\documentclass[reqno,12pt,a4paper]{amsart}
\usepackage{amsmath,amsthm,amsfonts,amssymb,amscd,amstext}
\usepackage[dvips]{graphicx}
\usepackage{psfrag,euscript,a4wide,stmaryrd,wasysym}

\usepackage[final]{epsfig}

\theoremstyle{plain}
\newtheorem{maintheorem}{Theorem}

\newtheorem{theorem}{Theorem}[section]
\newtheorem{proposition}[theorem]{Proposition}

\newtheorem{lemma}[theorem]{Lemma}
\newtheorem{claim}[theorem]{Claim}
\newtheorem{conjecture}{Conjecture}
\theoremstyle{definition}

\newtheorem{remark}[theorem]{Remark}
\newtheorem{definition}[theorem]{Definition}

\newcommand{\RR}{{\mathbb R}}
\newcommand{\NN}{{\mathbb N}}
\newcommand{\ZZ}{{\mathbb Z}}
\newcommand{\DD}{{\mathbb D}}
\newcommand{\sS}{{\mathbb S}}
\newcommand{\TT}{{\mathbb T}}
\newcommand{\CC}{{\mathbb C}}

\newcommand{\vfi}{\varphi}

\newcommand{\vr}{\varphi}

\newcommand{\diam}{\operatorname{diam}}
\renewcommand{\epsilon}{\varepsilon}
\renewcommand{\rho}{\varrho}
\newcommand{\dist}{\operatorname{dist}}

\newcommand{\Leb}{\operatorname{Leb}}
\newcommand{\var}{\operatorname{var}}

\newcommand{\const}{\operatorname{const}}

\newcommand{\cN}{\EuScript{N}}
\newcommand{\cC}{\EuScript{C}}
\newcommand{\cT}{\EuScript{T}}
\newcommand{\cF}{\EuScript{F}}

\newcommand{\cS}{\EuScript{S}}
\newcommand{\cI}{\EuScript{I}}

\newcommand{\U}{\EuScript{U}}
\newcommand{\cP}{\EuScript{P}}
\newcommand{\cO}{\EuScript{O}}

\newcommand{\cL}{\EuScript{L}}
\newcommand{\bT}{\mathbb{T}}

\providecommand{\abs}[1]{\lvert#1\rvert}

\providecommand{\Norm}[1]{\interleave#1\interleave}

\def \fX {{\mathfrak X}}

\def \fD {{\mathfrak D}}

\title[Multidimensional Rovella-like attractors]
{Multidimensional Rovella-like attractors}

\author{V. Araujo, A. Castro, M. J. Pacifico and V. Pinheiro}

\date{\today}

\begin{thanks} {V.A. and M.J.P. were partially supported by
    CNPq-Brazil/FAPERJ-Brazil/Pronex Dynamical Systems. }
\end{thanks}
\begin{document}

\address{V\'\i tor Ara\'ujo,
  Instituto de Matem\'atica,
  Universidade Federal do Rio de Janeiro C. P. 68.530,
  21.945-970, Rio de Janeiro, RJ-Brazil}
\email{vitor.araujo@im.ufrj.br}

\address{Armando de Castro, Departamento de Matem\'atica, Universidade Federal da Bahia\\
  Av. Ademar de Barros s/n, 40170-110 Salvador, Brazil.}
\email{armandomat@yahoo.com.br}

\address{Maria Jos\'e Pacifico, Instituto de Matem\'atica,
Universidade Federal do Rio de Janeiro,
C. P. 68.530, 21.945-970 Rio de Janeiro, Brazil}
\email{pacifico@im.ufrj.br}

\address{Vilton  Pinheiro, Departamento de Matem\'atica, Universidade Federal da Bahia\\
Av. Ademar de Barros s/n, 40170-110 Salvador, Brazil.}
\email{viltonj@ufba.br}

\begin{abstract}
% The geometrical model for the flow generated by
% the Lorenz equations is based on a modification of a linear
% flow in $\RR^3$ presenting a singularity $\sigma$ with real
% eigenvalues $\lambda_2 <  \lambda_3 < 0 < - \lambda_3 < \lambda_1 $.
% When the construction is made changing the above relation to
% $\lambda_2 >  \lambda_3 > 0 > - \lambda_3 > \lambda_1 $ the resulting
% geometric model is said Rovella-like.
  We present a multidimensional flow exhibiting a
  Rovella-like attractor: a transitive invariant set with a
  non-Lorenz-like singularity accumulated by regular orbits
  and a multidimensional non-uniformly expanding invariant
  direction.  Moreover, this attractor has a physical
  measure with full support and persists along certain
  submanifolds of the space of vector fields.  As in the
  $3$-dimensional Rovella-like attractor, this example is
  not robust.
  As a sub-product of the construction we obtain a new class
  of multidimensional non-uniformly expanding endomorphisms
  without any uniformly expanding direction, which is
  interesting by itself.
  Our example is a suspension (with singularities) of this
  multidimensional endomorphism.
  % upon which we build up a singular flow having this map as
  % a first return map to a global cross-section.
% , whose
%   suspension provides a Rovella-like attractor, which is
%   partially hyperbolic, and whose quotient over stable
  % To prove non-uniform expansion for this class of
  % multidimensional endomorphisms with critical points, we
  % analyze the return of the orbits of the multidimensional
  % critical point to itself under assumptions of slow
  % recurrence and asymptotic expansion (akin to the work of
  % Benedicks-Carleson).
  % we apply arguments of shadowing by iterates of the
  % critical values (akin to the work of Benedicks-Carleson),
  % and expansion away from a neighborhood of the critical
  % set. In this way

  % and prove non-uniform expansion for this class of
  % endomorphisms.
 \end{abstract}


\renewcommand{\subjclassname}{\textup{2000} Mathematics Subject Classification}


\maketitle

\setcounter{tocdepth}{2}

\tableofcontents

%\newpage

\section{Introduction and statements of the results}
\label{sec:introd-statem-result}

Flows in compact (two-dimensional) surfaces are very well
understood since the groundbreaking work of
Peixoto~\cite{peixoto59,peixoto62}. A theory of
three-dimensional flows has steadily been developing since
the characterization of robust invariant sets in
\cite{MPP04}.  The study of singular attractors for flows in
higher dimensions is, however, mostly open; see
\cite{BPV97,MeMor06} and the recent book \cite{AraPac2010}
and references therein.  % Extensions of these techniques for
% higher-dimensional flows are well under way, see
% e.g. \cite{MeMor06}.

For non-robust but persistent attractors, like the Rovella
attractor presented in \cite{Ro93}, there is still no
higher-dimensional analogue. Other previous constructions
yielding similar behavior can be found in,
e.g., \cite{Smnv78,Smnv78a,LybmvZks,TurShlnk87,GmbdTress88,GmdPrcThmTres}.
Of course, one may embed the usual Rovella attractor (also
known as \emph{contracting Lorenz} attractor) into flows in
any dimension, just by “multiplying by a strong
contraction” (the attractor is contained in a
three-dimensional submanifold, which is invariant and
normally contracting for the flow).  But this procedure
leads to flows without new dynamical phenomena.  Moreover,
it remains an open problem, since the introduction of the
contracting Lorenz models about two decades ago, whether
persistent non-robust attractors for singular flows may
contain singularities with multidimensional expanding
directions.

Here we present a positive solution to this problem.

In a rather natural way, since low dimensional dynamics is
much better understood than the dynamics of general systems,
the techniques used in the mathematical analysis of
three-dimensional attractors for flows frequently depend on
a dimensional reduction through projection along a stable
manifold inside some Poincar\'e cross-section. This method
yields a one-dimensional system whose dynamics can be nearly
completely understood as well as the dynamics of its small
perturbations; see e.g. \cite{MartdMelo01,MS93}.

In this work we start a rigorous study of a proposed higher
dimensional analogue of the three-dimensional Rovella
attractor. Other examples of higher dimensional chaotic
attractors have been recently presented, e.g. by Bonatti,
Pumari\~no and Viana in \cite{BPV97} and by Shilnikov and
Turaev in \cite{ST98}, but these are robust sets, while our
construction leads to a persistent non-robust singular
attractor.

In \cite{BPV97} the authors define a uniformly expanding map
on a higher-dimensional torus, suspend it as a time-one map
of a flow, and then singularize the flow adding a
singularity in a convenient flow-box. This procedure creates
a new dynamics on the torus presenting a multidimensional
version of the one-dimensional expanding Lorenz-like map,
and a flow with \emph{robust multidimensional Lorenz-like
  attractors}: the singularity contained in the attractor
may have any number of expanding eigenvalues, and \emph{the
  attractor remains transitive in a whole neighbourhood of
  the initial flow.} The construction in \cite{ST98} is also
robust but yields a quasi-attractor: the attracting
invariant set is not transitive but it is the maximal chain
recurrence class in its neighborhood.
% \begin{figure}[htbp]
% \centering
% \includegraphics[width=8cm]{Lorenz4d.eps}
% \caption{A sketch of the construction of a
%   singular-attractor in higher dimensions}
% \label{fig:sing-attractor-4d}
% \end{figure}

For the case of the Lorenz attractor and singular-hyperbolic
attractors in general, see e.g.~\cite{viana2000i} and
\cite{AraPac2010}. In this class the equilibrium accumulated
by regular orbits inside the attractor has real eigenvalues
$\lambda_{ss}<\lambda_s<0<\lambda_u$ satisfying
$\lambda_{s}+\lambda_u>0$ and a type of global (collection
of) cross-section(s) can be defined endowed with invariant
stable laminations.
\begin{figure}[htbp]
\psfrag{S}{$\Sigma$}
  \includegraphics[height=5cm,width=8cm]{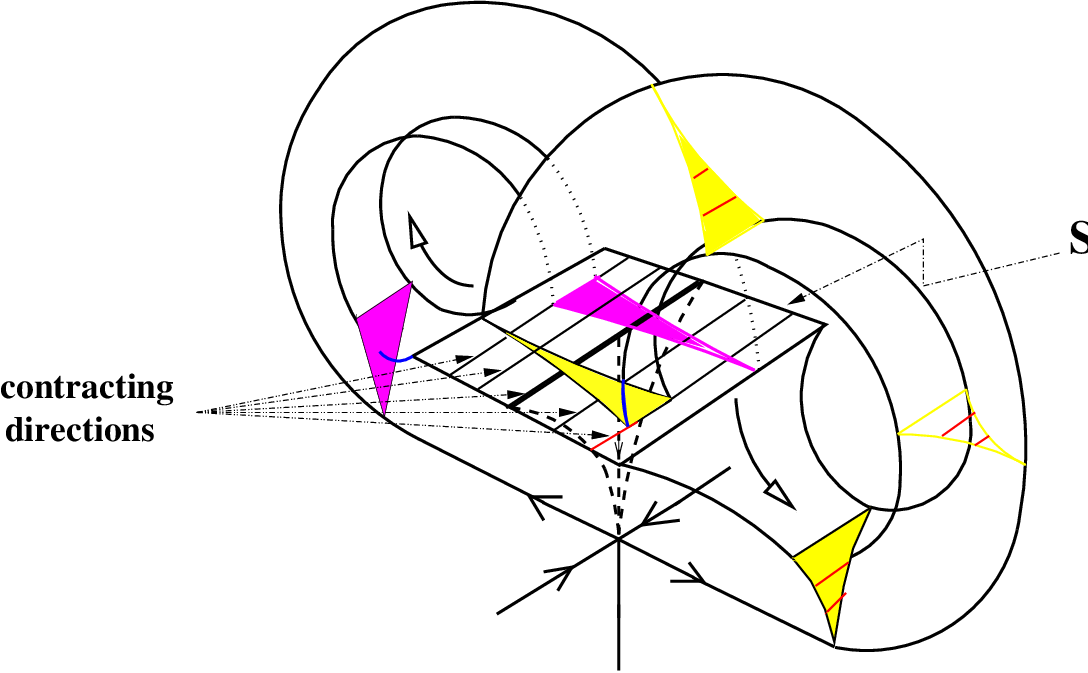}
  $\quad$
  \includegraphics[width=5cm]{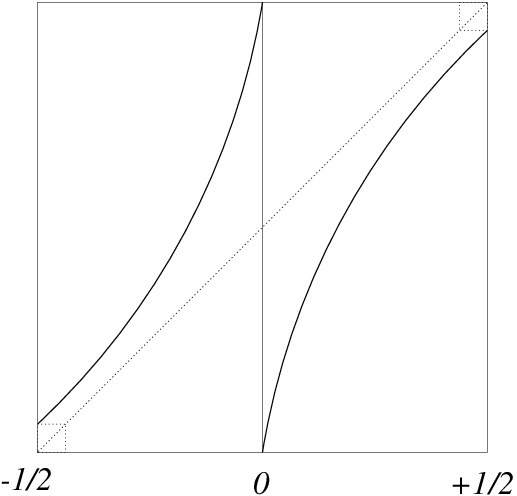}
  \caption{On the left, the geometric Lorenz attractor with
    the contracting directions on the cross-section
    $\Sigma$; on the right, the Lorenz one-dimensional
    transformation.}
  \label{fig:geom-lorenz}
\end{figure}
The quotient of the return map to the global cross-section
over the stable leaves (see Figure~\ref{fig:geom-lorenz})
is the one-dimensional Lorenz transformation.
% \begin{figure}[htpb]
% \includegraphics[width=4cm]{L1D.eps}
% \caption{\label{L1D}{The Lorenz one-dimensional transformation.}}
% \end{figure}
Our goal is to construct a flow such that ``after the
identification by the stable directions'', the first return
map in a certain cross section $M$ is a multidimensional
version of the one-dimensional map of the contracting Lorenz
model (or Rovella attractor); see~\cite{Ro93}.

\begin{figure}[htpb]
  \centering
  \includegraphics[width=14cm]{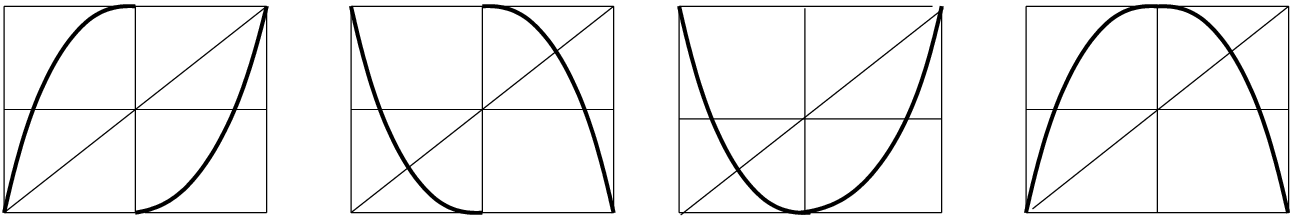}
  \caption{Several possible cases for the one-dimensional map for the
    contracting Lorenz model}
  \label{fig:contractingLorenz}
\end{figure}

A Rovella-like attractor is the maximal invariant set of a
geometric flow whose construction is very similar to the one
that gives the geometric Lorenz attractor,
\cite{GW79,ABS77,AraPac2010}, except for the fact that the
eigenvalues relation $\lambda_u + \lambda_s > 0$ there is
replaced by $\lambda_u + \lambda_s < 0$, where $\lambda_u>
0$ and $\lambda_s$ is the weakest contractive eigenvalue at
the origin. We remark that, unlike the one-dimensional
Lorenz map obtained from the usual construction of the
geometric Lorenz attractor, a one-dimensional map associated
to the contracting Lorenz model has a criticality at the
origin, caused by the eigenvalue relation $\lambda_u +
\lambda_s<0$ at the singularity.  In
Figure~\ref{fig:contractingLorenz} we present some possible
one-dimensional maps obtained through quotienting out the
stable direction of the return map to the global
cross-section of the geometric model of the contracting
Lorenz attractor, as in Figure~\ref{fig:geom-lorenz}.

The interplay between expansion away from the critical point
with visits near the criticality prevents this system to
have uniform expansion and also prevents robustness, that
is, the attractor is not transitive on a whole neighborhood
of the original flow.  However the Rovella attractors are
\emph{persistent}. 

We say that an attractor $\Lambda$ of a vector field
$X\in\fX^3(M)$ is \emph{$k$-dimensionally persistent}, if
there exists a $C^\infty$ submanifold $\cP$ of $\fX^3(M)$ with
codimension $k$ and $\Lambda$ admits a neighborhood $U$ in
$M$ such that for every $Y\in\cT:=\{Y\in\fX^3(M):
\cap_{t\ge0}Y^t(U) \text{ is transitive } \}$ we have
\begin{align*}
  \lim_{r\to0} \frac{m(\cT\cap S\cap B_r(Y) )}{m(S\cap B_r(Y))}=1
\end{align*}
where $S$ is any $k$-dimensional submanifold of $\fX^3(M)$
intersecting $\cP$ transversely, $m$ is $k$-dimensional
Lebesgue (volume) measure in $S$ and $B_r(Y)$ is the ball of
radius $r>0$ in $S$ in the $C^3$ topology.

Rovella in~\cite{Ro93} showed that this class of
three-dimensional attractors is $2$-dimensionally persistent
in the $C^3$ topology.  That is, \emph{for generic
  parameterized families of vector fields passing through
  the original vector field $X_0$, the parameters
  corresponding to transitive attractors, with the same
  features of the Rovella attractor, form a positive
  Lebesgue measure subset and full density at $X_0$.} We
stress that Rovella-like and Lorenz-like attractors are
rather natural dynamical models since they appear in the
generic unfolding of resonant double homoclinic loops; see
e.g. \cite{Ro89,Ro92,Ro2000,MPS05,MPsM}.

In this work we provide a multidimensional counterpart of
this result. First we obtain the attractor as follows.

\begin{maintheorem}\label{mthm:pri}
  For any dimension $m=k+5$, $k \in\ZZ^+$, there exist a
  $C^{\infty}$ vector field $X \in {\fX}^\infty(M^m)$ on a
  $m$-dimensional manifold exhibiting a singular-attractor
  $\Lambda$, containing a pair of hyperbolic singularities
  $s_0,s_1$ \emph{with different indexes} in a trapping
  region $U$. Moreover
  \begin{enumerate}
  \item there exists a map
    $R:\Sigma\circlearrowleft$ on a $(k+4)$-dimensional
    cross-section $\Sigma$ of the flow of $X$ such that
    \begin{enumerate}
    \item the set $\Lambda$ is the suspension of an attractor
      $\Lambda_\Sigma\subset\Sigma$ with respect to $R$;
    \item $R$ admits a $3$-dimensional stable direction $E^s$
      and $(k+1)$-dimensional center-unstable direction $E^c$
      such that $E^s\oplus E^c$ is a partially hyperbolic
      splitting of $T\Sigma\mid\Lambda_\Sigma$;
    \item $\Lambda_\Sigma$ supports a physical measure $\nu$
      for $R$.
    \end{enumerate}
  \item $\Lambda$ is the
    support of a physical hyperbolic measure $\mu$ for the
    flow $X^t$: the ergodic basin of $\mu$ is a positive
    Lebesgue measure subset of $U$ and every Lyapunov
    exponent of $\mu$ along the suspension of the bundle
    $E^c$ is positive, except along the flow direction.
  \end{enumerate}
\end{maintheorem}

To prove Theorem~\ref{mthm:pri} we follow the same strategy
of \cite{BPV97} with two main differences.

On the one hand, since we aim at a multidimensional version
of the one-dimensional map of the contracting Lorenz model,
we have to deal with critical regions, that is, regions
where the derivative of the return map to a global
cross-section vanishes.  Because of this, proving the
existence of non-trivial attractors for the flow arising
from such construction requires a more careful analysis.
Indeed, as in the one-dimensional case, depending on the
dynamics of the critical region, every attractor for the
return map may be periodic (trivial).

On the other hand, our construction leads to the presence of
a pair of hyperbolic saddle equilibria accumulated by
regular orbits inside the attractor but \emph{with different
  indexes} (the dimension of their stable manifolds). This
feature creates extra difficulties for the analysis of the
possible dynamics arising under small perturbations of the
flow.

Typically, when the critical region is non-recurrent (which
corresponds to Misiurewicz maps in one-dimensional
dynamics), most of the difficulties introduced by the
critical region can be bypassed. That is one of the main
reasons for us to construct a kind of multidimensional
Misiurewicz dynamics.  In general, such critical regions in
dimension greater than one are sub-manifolds, and one cannot
rule out that they intersect each other under the action of
the dynamics.  Albeit this, we shall exhibit a class of
multidimensional Misiurewicz-like endomorphisms that appears
naturally in a flow dynamics; see Theorem~\ref{thm:seg} in
Section~\ref{sec:theconstruction}. This is an example of
non-uniformly expading dynamics in higher dimensions which
\emph{does not exhibit any uniformly expanding direction}
and is conjugated to a skew-product of a quadratic map with
an expanding map, with the exception of at most two orbits.

This is the basic dynamics which we modify to obtain the
return map $R$ to a cross-section of the flow, exhibiting an
expanding invariant torus $\bT_1^k$ that will absorve the
image of the critical region after the singularization of
the associated flow.  \emph{By topological reasons, this map
  can not be seen as a time-one map of a suspension flow:
  locally its degree is not constant.}  To bypass this new
difficulty, we realize this map as a first return map of a
flow \emph{with singularities} (after identification by
stable directions).  Afterwards, we singularize a periodic
orbit of this flow, i.e., we introduce a new singularity $s$
of saddle-type, with $(k+1)$-dimensional unstable manifold
and $4$-dimensional stable manifold. Moreover, all the
eigenvalues of $s$ are real and if $\sigma_{s,i}$, and
$\sigma_{u,j}$ denote the stable and the unstable
eigenvalues at $s$ respectively, then $\max\{\sigma_{s,i}\}
+ \max\{\sigma_{u,j}\} < 0$ for $1\le i \le 4$ and $1 \le j
\le k+1$.  We say that this kind of singularity is a {\emph
  {Rovella-like singularity}}.  The resulting flow will
present a multidimensional transitive \emph{Rovella-like
  attractor}, supporting a physical measure, as stated in
Theorem~\ref{mthm:pri}.

The existence of the physical/SRB measure for the original
flow is obtained taking advantage of the fact that, through
identification of stable leaves, we can project the dynamics
of the first return map $R$ of the flow to a global
cross-section into a one-dimensional transformation with a
Misiurewicz critical point.

We point out that \emph{the analysis of the dynamics of most
  perturbations of our flow cannot be easily reduced
  (perhaps not at all) to a one-dimensional model}. This
indicates that intrinsic multidimensional tools should be
developed to fully understand this class of flows. Thus
extra difficulties arise to verify that this kind of
multidimensional contracting Lorenz attractor is persistent.
We obtain the following partial result in this direction.
\begin{maintheorem}
  \label{mthm:codimension}
  For any $k\in\ZZ^+$, there exists a $(k+2)$-codimension
  submanifold $\cP$ of the space of $C^2$ vector fields
  $\fX^2(M)$ such that
  \begin{enumerate}
  \item 
    the vector field $X$ from Theorem~\ref{mthm:pri} belongs
    to $\cP$ and, for all $Y\in\cP$ in a neighborhood $\U$
    of $X$ in $\fX^2(M)$, the Poincar\'e return map $R_Y$ to
    the cross-section $\Sigma$ admits a $3$-dimensional
    strongly contracting $C^\gamma$ foliation $\cF$, for
    some $\gamma>1$.  The induced map $g_Y$ on the quotient
    of $\Sigma$ over $\cF$ is a $C^{\gamma}$ endomorphism on
    a cylinder $[-1,1]\times\TT^k$.
  \item 
    for vector fields $Y\in\U\setminus\cP$, the
    Poincar\'e return map $R_Y$ to the cross-section
    $\Sigma$ admits a one-dimensional strongly contracting
    $C^\gamma$ foliation $\cF$, for some $\gamma>1$.  The
    induced map $g_Y$ on the quotient of $\Sigma$ over $\cF$
    is a $C^{\gamma}$ endomorphism on a manifold
    diffeomorphic to the unit ball in $\RR^{k+3}$.
  \item for a vector field $Y\in\U\cap\cP$, %  where $\U$ is
    % the neighborhood of $X$ and $\cP$ the submanifold in
    % $\fX^2(M)$ introduced at Theorem~\ref{mthm:codimension},
    \emph{if} the quotient map $g_Y$ sends the critical set
    inside the stable manifold of the sink $p(Y)$,
    \emph{then} this stable manifold contains the trapping
    region $U$ except for a zero Lebesgue measure subset.
  \end{enumerate}
\end{maintheorem}

% Considering the perturbation of this flow along parametrized
% families, we can show the existence of many parameters for
% which nearby flows exhibit an attractor with a unique
% physical measure. For this we use the construction of a
% partial hyperbolic Poincar\'e return map to a cross-section
% of the flow and take advantage of this weak hyperbolicity
% through a multidimensional extension of arguments of
% Benedicks-Carleson type. That is, we analyze the return of
% the orbits of the multidimensional critical set to itself
% under certain natural assumptions of slow recurrence and
% asymptotic expansion.

Next we give the precise definitions and concepts involved
in the previous statements.

%%%%%%%%%%%%%%%%%%%%%%%%%%%%%%%%%%%%%%%%%%%%%%%%%%%%%%%

\subsection{Preliminary definitions and conjectures}
\label{sec:preliminaryresults}

In what follows $M$ is a compact boundaryless finite
dimensional manifold and ${\fX}^1(M)$ is the set of $C^1$
vector fields on $M$, endowed with the $C^1$ topology.  From
now on we fix some smooth Riemannian structure on $M$ and an
induced normalized volume form $m$ that we call Lebesgue
measure.  We write also $\dist$ for the induced distance on
$M$ and $\|\cdot\|$ for the induced Riemannian norm on $TM$.
Given $X \in {\fX}^1(M)$, we denote by $X^t$, $t \in \RR$
the flow induced by $X$, and if $x \in M$ and $[a,b]\subset
\RR$ then $X^{[a,b]}(x)= \{X^t(x), a\leq t \leq b\}$.

We say that a differentiable map $f:M\circlearrowleft$ is
$C^{1+}$ if the derivative $Df$ of $f$ is \emph{H\"older}:
there are $\alpha,C>0$ such that for every $x\in M$ we can
find parametrized neighborhoods $U=\phi(U_0)$ of $x$ and
$V=\psi(V_0)$ of $f(x)$ in $M$, where $U_0,V_0$ are
neighborhoods of $0$ in $\RR^{\dim(M)}$ and $\phi,\psi$ are
parametrizations of $M$, such that for all $y_1,y_2\in U$
\begin{align*}
  \|D(\psi^{-1}\circ f\circ \phi)(y_1)-D(\psi^{-1}\circ
  f\circ \phi)(y_2)\|
  \le
  C \dist(y_1,y_2)^\alpha.
\end{align*}

A point $p\in M$ is a \emph{periodic point} for $X^t$ if
$X(p)\neq0$ and there exists $\tau>0$ such that
$X^\tau(p)=p$. The minimal value of $\tau$ such that
$X^\tau(p)=p$ is the \emph{period} of $p$. If $p$ is a
periodic point, we also say that the orbit $\cO(p)=\{
X^t(p): t\in\RR\}$ of $p$ is a \emph{periodic orbit}. A
\emph{singularity} $\sigma$ is an equilibrium point of
$X^t$, that is, $X(\sigma)=0$. If $X(p)\neq0$ then we say
that $p$ is a \emph{regular point} and its orbit $\cO(p)$ is
a regular orbit.

Let $\Lambda$ be a compact invariant set of $X\in
{\fX}^1(M)$.  We say that $\Lambda$ is an \emph{attracting
  set} if there exists an \emph{trapping region}, i.e. an
open set $U\supset \Lambda$ such that
$\overline{X^t(U)}\subset U$ for $t>0$ and $ \Lambda
=\bigcap_{t\in \RR}X^t(U)$. Here $\overline{A}$ means the
topological closure of the set $A$ in the manifold we are
considering.

We say that an attracting set
$\Lambda$ is \emph{transitive} if it is equal to the
$\omega$-limit set of a regular $X$-orbit. We recall that
the $\omega$-limit set of a given point $x$ with respect to
the flow $X^t$ of $X$ is the set $\omega(x)$ of accumulation
points of $(X^t(x))_{t>0}$ when $t\to+\infty$.  An
\emph{attractor} is a transitive attracting set and a
\emph{singular-attractor} is an attractor which contains
some equilibrium point of the flow.  An attractor is
\emph{proper} if it is not the whole manifold.  An invariant
set of $X$ is \emph{non-trivial} if it is neither a periodic
orbit nor a singularity.

\begin{definition}
\label{d.dominado}
Let $\Lambda$ be a compact invariant set of a $C^{1+}$ map
$f:M\circlearrowleft$ , $c>0$, and $0 < \lambda < 1$.  We
say that $\Lambda$ has a $(c,\lambda)$-dominated splitting
if the bundle over $\Lambda$ splits as a $Df$-invariant sum
of sub-bundles $ T_\Lambda M=E^s\oplus E^{cu}, $ such that
for all $n\in\ZZ^+$ and each $x \in \Lambda$
\begin{equation}\label{eq.domination}
\|Df^n \mid E^s_x\| \cdot \|(Df^{n} \mid E^{cu}_{x})^{-1}\| < c \,
\lambda^n.
\end{equation}
\end{definition}

% The domination condition \eqref{eq.domination} implies that
% the flow direction is contained in the $E^{cu}$ sub-bundle.

We say that a $f$-invariant subset $\Lambda$ of $M$ is
\emph{partially hyperbolic} if it has a
$(c,\lambda)$-dominated splitting, for some $c>0$ and
$\lambda\in(0,1)$, such that the sub-bundle $E^s$ is
uniformly contracting: for all $n\in\ZZ^+$ and each $x \in
\Lambda$ we have
\begin{align}\label{eq:unif-contr}
  \|Df^n \mid E^s_x\| < c \, \lambda^n.
\end{align}

We denote by $\overline{A}$ the topological closure of the
set $A\subset M$ in what follows.

\begin{definition}
  \label{def:NUE}
  A $C^{1+}$ map $g:M\to M$ is non-uniformly expanding if
  there exists a constant $c>0$ such that
  \begin{align*}
    \liminf_{n\to+\infty}
    \frac1n\sum_{j=0}^{n-1}\log\|(Dg(g^j(x))^{-1}\| \le -c<0 \quad\text{for
      Lebesgue almost every}\quad x\in M.
  \end{align*}
\end{definition}

Non-uniform expansion ensures the existence of absolutely
continuous invariant probability measures under some mild
extra assumptions on $g$, see below. One property of such
measures is that they have a ``large ergodic basin''.

\begin{definition}{(physical measure and Ergodic basin for flows.)}
  We say that an invariant probability measure $\mu$ is a
  {\em physical measure} for:
  \begin{itemize}\item 
    the flow given by a field $X:M \to TM$ if there exists a
    positive Lebesgue measure $B(\mu)\subset M$ such that
    for all $x \in B(\mu)$ we have
$
\frac{1}{T} \int_{0}^T \vr \circ X^t(x) dt
\xrightarrow[T\to+\infty]{} \int_M \vr d\mu
$
for all continuous functions $\vfi:M\to\RR$;
\item  the map $g:M\to M$ if there exists a positive
  Lebesgue measure $B(\mu)\subset M$ such that for each $x
  \in B(\mu)$ we have
$
\frac{1}{n} \sum_{j=0}^{n-1} \vfi(g^j(x))
\xrightarrow[n\to+\infty]{} \int_M \vfi\, d\mu
$
for all continuous functions $\vfi:M\to\RR$.
\end{itemize}
The set $B(\mu)$ is
called the {\em ergodic basin} of $\mu$.
\end{definition}

Note that Theorems~\ref{mthm:pri} and~\ref{mthm:codimension}
show that the singular-attractor $\Lambda$ is not robustly
transitive not even robust: there exist arbitrary small
perturbations $Y$ of the vector field $X$ such that the
orbits of the flow of $Y$, in a full Lebesgue measure set of
points in $U$, converge to a periodic attractor (a periodic
sink for the flow). The singular-attractor $\Lambda$ is not
partially hyperbolic in the usual sense
\eqref{eq.domination} and~\eqref{eq:unif-contr} adapted to
the flow setting, since we can only define the splitting on
the points of $\Lambda$ which do not converge to the
singular points $s_0,s_1$, that is, we exclude the stable
set of the singularities within $\Lambda$. The remaining
set, however, has full measure with respect to $\mu$. Since
the equilibria $s_0,s_1$ contained in $\Lambda$ are
hyperbolic with different indexes (i.e. the dimension of
their stable manifolds) we believe the following can be
proved.

\begin{conjecture}\label{conj:noextension}
  It is not possible to extend the dominated splitting on
  $\Lambda$ away from equilibria to the equilibria $s_0,s_1$
  which belong to $\Lambda$.
\end{conjecture}

In addition, the attractor $\Lambda$ for $X$ in $U$ is such
that the Jacobian along any $2$-plane $P_x$, inside the
central subbundle $E^c_x$ for Lebesgue almost all $x\in U$,
is \emph{asymptotically expanded but not uniformly
  expanded}, that is we can find a constant $c>0$ such that
$
  \lim_{t\to+\infty}\frac1t\log|\det DX^t\mid P_x|\ge c,
$
but $\frac1t\log|\det DX^t\mid P_x|$ can take an arbitrary
long time (depending on $x$) to become positive.

For the remaining cases, not analyzed in
Theorem~\ref{mthm:codimension}, we conjecture that the
quotient map behaves in a similar way to the typical
perturbations of a smooth one-dimensional multimodal map.

\begin{conjecture}
  \label{conj:physicalmeas}
  For a vector field $Y\in\U$, where $\U$ is the
  neighborhood of $X$ in $\fX^2(M)$ introduced in
  Theorem~\ref{mthm:codimension}, if the quotient map $g_Y$
  \emph{does not send} the critical set inside the stable
  manifold of the sink $p(Y)$, then the complement of this
  stable manifold contains the basin of a physical measure
  for $g_Y$ whose Lyapunov exponents are positive.
\end{conjecture}

To make any progress in the understanding of this
conjecture one needs to study the interplay between the
critical set, the expanding behavior in some regions of the
space, and the stable manifold of the sink, in a higher
dimensional setting. We believe this will demand the
development of new ideas in dynamics and ergodic theory.

%%%%%%%%%%%%%%%%%%%%%%%%%%%%%%%%%%%%%%%%%%%%%%%%%%%%%%%%%

\subsection{Organization of the text}
\label{sec:open-questi}

We present the construction of the vector field $X$ in
stages in Section~\ref{sec:theconstruction}. We start by
constructing a non-uniformly expanding higher dimensional
endomorphism with critical points in
Section~\ref{sec:an-example-nue}. This provides a starting
point for the Poincar\'e return maps
$R:\Sigma\circlearrowleft$ of the statement of
Theorem~\ref{mthm:pri}. Then we adapt this first
construction to become the quotient of the return map of the
flow we will construct, in
Section~\ref{sec:unpert-basic-dynamic}.  In
Section~\ref{sec:unpert-singul-flow} we start the
construction of the singular flow we are interested in. This
is done again in stages, and here we obtain a first
candidate.  Next we obtain the vector field $X$ after
perturbing the candidate in Section~\ref{sec:pert_flow}.

We study the properties of $X$ and its unfolding in
Section~\ref{sec:conseq-unfold-x}. We prove the existence of
the dominated splitting for the return map to the
cross-section and describe the construction of the physical
measure with positive multidimensional Lyapunov exponents in
Section~\ref{sec:x-chaotic-with}, completing the proof of
Theorem~\ref{mthm:pri}. The details on the existence of the
physical measure for the return map are left for
Section~\ref{sec:higher-dimens-misiur}, where
``Benedicks-Carleson type'' arguments are adapted to our
higher-dimensional setting.

The unfolding of the vector field $X$ is studied in
Section~\ref{sec:unfold-x}, where items (1-2) of
Theorem~\ref{mthm:codimension} are proved and the argument for the
proof of item (3) is described. The proof of this
last item is given in Section~\ref{sec:full-basin-attract}.

We end with two
appendixes. Appendix~\ref{sec:hyperb-neighb-constr} provides
an adaptation of a major technical tool to our setting to
prove the existence of an absolutely continuous measure for
higher dimensional non-uniformly expanding maps.
Appendix~\ref{sec:isotopy} provides a topological argument
for the existence of a certain isotopy we need during the
construction of the vector field $X$.

\subsection*{Acknowledgments}

This paper started among discussions with A.C.J. and
V.P. during a (southern hemisphere) summer visit V.A. and
M.J.P. payed to Universidade Federal da Bahia, at
Salvador. We thank the hospitality and the relaxed and
inspiring atmosphere of Bahia.

%%%%%%%%%%%%%%%%%%%%%%%%%%%%%%%%%%%%%%%%%%%%%%%%%%%%%%%%%%%%

\section{The construction}
\label{sec:theconstruction}

Here we prove Theorem~\ref{mthm:pri}. We start by 
%proving
%the following result.
%
%\begin{theorem} \label{thm:seg} For any $k\in\ZZ^+$, there
 % exist a non-uniformly expanding map $g$ of class $C^{1+}$
 % on a $(k+1)$-manifold $N$, without any uniformly expanding
  %directions, admitting an absolutely continuous invariant
  %probability measure with full ergodic basin in $N$, whose
  %Lyapunov exponents are all positive.
%\end{theorem}
%
%This provides 
providing 
an example of non-uniformly expading dynamics
in higher dimensions which is interesting by itself since it
\emph{does not exhibit any uniformly expanding direction}. 
%  This will be used to
% construct the attractor for the flow of
% Theorem~\ref{mthm:pri}.
Then we adapt this example to obtain the Poincar\'e return map
$R$ in the statement of Theorem~\ref{mthm:pri}.
At the end of this section, we show how this yields a
construction of a singular attractor with $k+1$ non-uniformly
expanding directions on a compact boundaryless
$(k+5)$-dimensional manifold $M$.

%%%%%%%%%%%%%%%%%%%%%%%%%%%%%%%%%%%%%%%%%%%%%%%%%%%%%%%

\subsection{An example of non-uniformly expanding dynamics
  in high dimension}
\label{sec:an-example-nue}

We start by defining a non-uniformly expanding endomorphism
of a $(k+1)$- dimensional manifold $N$.

\begin{figure}[htpb]
  \centering\psfrag{s}{$-1$}\psfrag{n}{$+1$}\psfrag{t}{$\TT^k$}
  \psfrag{N}{$N$}
  \includegraphics[width=6cm]{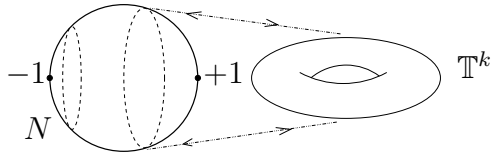}
  \caption{The "torusphere": each parallel is a $k$-torus.}
  \label{fig:torusph-each-meridi}
\end{figure}

Let $\Upsilon:\RR\times \CC^{k} \to \RR\times \CC^{k}$ be
given by $(t, z)\mapsto (\cos(\pi t), z \sin(\pi t))$.  We
consider $\TT^k=\sS^1 \times \overset{k}{\dots} \times
\sS^1$ where $\sS^1=\{z\in\CC: |z|=1\}$ and let $N=
\Upsilon([-1,1]\times \TT^{k})$.  Clearly $N$ has a natural
differential structure: $N=G^{-1}(\{1,\dots,1\})$ for
$$
G:\RR\times\CC^{k}\to\RR^k:
(t,z_1,\dots,z_k)=(t^2+|z_1|^2,\dots,t^2+|z_k|^2)
$$
which we name "torusphere" and is a manifold of dimension
$k+1$, see Figure~\ref{fig:torusph-each-meridi}.

We remark that this manifold is the boundary of
$M:=G^{-1}([0,1]^k)$, which is a ``solid torusphere'', that
is $M\cap(\{t\}\times\RR^k) \simeq \TT^k\times\DD$ for all
$-1<t<1$, where $\DD$ is the unit disk in $\CC$.
In what follows we write $I= [-1, 1]$.

Let  $g_0: I \circlearrowleft$ be a $C^{1+}$ non-flat
unimodal map with the critical points $c_0=0$ and
$c_1=g_0(c_0)=1$ as follows
$$
g_0(x)=
\begin{cases}
  \varsigma^+\left( 
    \left|\frac{x}2\left(1-\frac{x}2\right)\right|^\alpha
  \right) 
  & \text{if   } x\in[0,1]\\
  \varsigma^-(|x|^\alpha) & \text{if } x\in[-1,0)
\end{cases};
$$
for $C^\infty$ diffeomorphisms $\varsigma^\pm:[0,1]\to I$ such
that both $\varsigma^+$ and $\varsigma^-$ are monotonous
decreasing. Moreover we assume that \emph{the critical order
$\alpha$ is strictly between $1$ and $2$}, $1<\alpha<2$, and
that $g_0$ satisfies (see Figure~\ref{fig:1dnue}):
\begin{figure}[htbp]
\begin{center}
\psfrag{s}{$-1=p_0$}\psfrag{q}{$p_1$}\psfrag{n}{$+1$}
\psfrag{0}{$0$}
  \includegraphics[scale=0.3]{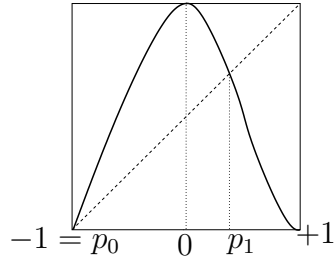}
\end{center}
\caption{The one-dimensional map $g_0$.}
  \label{fig:1dnue}
\end{figure}
\begin{enumerate}
\item $g_0(\pm1)=-1$ and $g_0$ has exactly two fixed points,
  namely, $p_0=-1$ and $0<p_1<1$;
\item $g_0'( p_0) > 1$ and $g_0'(p_1)<-1$.
\end{enumerate}

% Rescaling the interval $[0, \pi]$ to $I$ by means of $q(\theta):=
% (\pi- 2 \theta)/\pi$, we define $\hat g_0:= q^{-1} \circ g_0 \circ
% q$.
% Note that $|g'_0|= |\hat g'_0|$.

It is well known that these maps are non-uniformly
expanding and conjugated to the tent map (see the proof of
Lemma~\ref{le:alamane} in what follows). Thus, in particular,
they are topologically transitive and admit a unique
absolutely continuous probability measure supported on the
entire interval.

Now we define $f_1: \CC^k\to \CC^k$ by
$(w_1,\dots,w_k)\mapsto (w_1^2,\dots,w_k^2).$ Then the map
$g_0\times f_1: [-1, 1]\times \TT^k\circlearrowleft $
induces a $C^{1+\alpha}$ map $g : N \circlearrowleft$ by
$
g:=\Upsilon \circ (\hat g_0\times f_1 ) \circ \Upsilon^{-1}.
$
It is easy to see that $g$ takes meridians onto meridians
and parallels onto parallels of the torusphere. This means
that the derivative of $g$ preserves the directions
associated to such foliations.

Next we state and prove the main properties of the map $g$
defined above.  This is the main result of this subsection.

\begin{theorem} \label{thm:seg} For any $k\in\ZZ^+$, $g$ is
  a non-uniformly expanding map of class $C^{1+}$ on the
  $(k+1)$-manifold $N$, without any uniformly expanding
  directions, admitting an absolutely continuous invariant
  probability measure with full ergodic basin in $N$, whose
  Lyapunov exponents are all positive.
\end{theorem}

For the proof of theorem above we start with the following
\begin{claim}
  \label{claim:NUE}
  The map $g$ is a non-uniformly expanding map.
\end{claim}

Using polar coordinates, we can take a parametrization $h:
(-1,1) \times (0, 2\pi)^k \to N$ given by $ h(t,\Theta)=
h(t, \theta_1, \dots, \theta_k):= (\sin(\frac\pi2 t),
z(\Theta)\cos(\frac\pi2 t)), $ where $\Theta:= (\theta_1,
\dots, \theta_k)$ and
$$
z(\Theta):=
(\cos(\theta_1), \sin(\theta_1), \dots, \cos(\theta_k),
\sin(\theta_k)) \in \RR^{2k}.
$$
% $
% \cos(\theta_1) \cos(\frac\pi2 t),
% \sin(\theta_1) \cos(\frac\pi2 t),  \dots,\cos(\theta_k) \cos(\frac\pi2 t),
% \sin(\theta_k) \cos(\frac\pi2 t)),
% $$
The image of $h$ covers $N$ except for a null Lebesgue
measure set. % Let us write $\Theta:= (\theta_1, \dots,
% \theta_k)$, $z(\theta):= (\cos(\theta_1), \sin(\theta_1),
% \dots, \cos(\theta_k), \sin(\theta_k)) \in
% \RR^{2k}$.
Therefore, the expression of $g$ in these coordinates is
$$
\underbrace{\Big(\sin(\frac\pi2 t), \cos(\frac\pi2 t) z(\Theta)\Big)}_{:= x}
\longmapsto \Big(\sin( \frac\pi2 g_0(t)), \cos(\frac\pi2
g_0(t)) z(2\Theta) \Big).
$$
In the meridian directions, this implies that the derivative
$Dg(x): T_x N \to T_{g(x)} N$ takes the vector $v:=
\big(\cos(\frac\pi2 t), \sin(\frac\pi2 t) z(\Theta)\big)$ 
to 
$Dg(x) \cdot v=
g'_0(t) \cdot \big(\cos(\frac\pi2 g_0(t)), \sin(\frac\pi2 g_0(t))
z(2\Theta)\big)$. 

We adopt in $\RR \times \RR^{2k}$ (where $N$ is embedded)
the norm $\Norm{(t, z)}:= \sqrt{\abs{t}^2+ \|z\|^2/k}$,
where $\|z\|$ is the standard Euclidean norm in
$\RR^{2k}$. Therefore
\begin{align*}
  \frac{\Norm{Dg(x) \cdot v}}{\Norm{v}}
&= | g'_0(t)| \cdot
  \sqrt{\frac{\cos^2(\frac\pi2 g_0(t))+ \sin^2(\frac\pi2 g_0(t)) \cdot
      \|z(2\Theta)\|^2/k}{\cos^2(\frac\pi2 t)+ \sin^2(\frac\pi2 t) \cdot
      \|z(\Theta)\|^2/k} }
\\
&=
| g'_0(t)| \cdot \sqrt{\frac{\cos^2(\frac\pi2 g_0(t))+
    \sin^2(\frac\pi2 g_0(t))}{\cos^2(\frac\pi2 t)+
    \sin^2(\frac\pi2 t) } }= |g'_0(t)|.
\end{align*}
Along the directions of the parallels, given any $j\in \{1,
\dots, k\}$ the derivative $Dg(x)$ takes the vector $v_j:= (0,
\dots, -\cos(\frac\pi2 t)\sin(\theta_j), \cos(\frac\pi2 t)\cos(\theta_j), 0,
\dots 0)$ to the vector
$$
Dg(x) \cdot v_j= (0,  \dots, -2 \cos(\frac\pi2 g_0(t))\sin(2 \theta_j),
2\cos(\frac\pi2 g_0(t))\cos(2 \theta_j), 0, \dots 0).
$$
Therefore
$$
\frac{\||Dg(x) \cdot v_j|\|}{\||v_j|\|}
= 
2 \frac{|\cos(\frac\pi2 g_0(t))|
\cdot \sqrt{\sin^2(2\theta_j)+ \cos^2(2\theta_j)}}
{|\cos(\frac\pi2 t)|\cdot \sqrt{\sin^2(\theta_j)+
\cos^2(\theta_j)}} = 
2 \Big|\frac{\cos(\frac\pi2 g_0(t))}{\cos(\frac\pi2 t)}\Big|.
$$
We note that for $w= \sum_{j= 1}^k \alpha_j v_j$, the
relation $ \frac{\Norm{Dg(x) \cdot w}}{\Norm{w}}= 2
\Big|\frac{\cos(\frac\pi2 g_0(t))}{\cos(\frac\pi2 t)}\Big| $
also holds. Let us call $E_x$ the space generated by the
directions of the parallels through $x$ in $T_x N$, and
$m(x):= \inf_{w \in E_x} \frac{\Norm{Dg(x) \cdot
    w|}}{\Norm{w}}$ the minimum norm (or conorm) of $Dg(x)$
restricted to $E_x$. We have
\begin{align*}
  \frac 1 n \sum_{j= 0}^{n-1} \log m(g^j(x))
  &=
  \log(2) +
  \frac 1 n \sum_{j= 0}^{n-1} \log\Big|\frac{\cos(\frac\pi2
    g_0^{j+1}(t))}{\cos(\frac\pi2 g_0^{j}(t))}\Big|
= \log(2) +
  \frac1n
  \log\left|
      \frac{\cos(\frac\pi2 g_0^n(t))}{\cos(\frac\pi2 t)}
    \right|.
\end{align*}
Thus we will have $n^{-1} \sum_{j= 0}^{n-1} \log
m(g^j(x)) \ge \log 2$ whenever $|\cos(\frac\pi2 g_0^n(t))|\ge
|\cos(\frac\pi2 t)|$, which is true if, and only if,
$|g_0^n(t)|\le |t|$.

To conclude the proof of Claim~\ref{claim:NUE} we use
Lemma~\ref{le:alamane} below, whose proof we postpone to the
end of this subsection.
\begin{lemma}
  \label{le:alamane}
  Given a neighborhood $U$ of $c_0$, Lebesgue almost every
  orbit visits $U$ infinitely often.
\end{lemma}

Since $c_0=0$, Lemma~\ref{le:alamane} ensures that for every
given $t\in I\setminus\{0\}$
% $g_0(c_0)=c_1=1$, the continuity of $g_0$ implies that
% for Lebesgue almost every $t\in[-1,1]$
the inequality
$|g_0^n(t)|\le|t|$ is true for infinitely many values of
$n\ge1$. This implies that
\begin{align}
  \label{eq:NUEinf0}
  \liminf_{n\to+\infty}
  \frac1n\sum_{j=0}^{n-1}\log\Norm{(Dg\mid
    E_{g^j(x)})^{-1}} \le -\log 2<0 \quad\text{for Lebesgue
    almost every}\quad x.
\end{align}
Denoting $F_x$ the direction of the meridian at $T_x N$, for
$x=(t,\Theta)$, we also showed that
\begin{align}
  \label{eq:NUEinf1}
  \liminf_{n\to+\infty}
  \frac1n\sum_{j=0}^{n-1}\log\Norm{(Dg\mid F_{g^j(x)})^{-1}}
  =
  -\limsup\frac1n\sum_{j=0}^{n-1}\log|g_0^\prime(g_0^j(t))|<0
\end{align}
for Lebesgue almost every $t\in I$, which is strictly
negative by known results on unimodal maps (see
e.g. \cite{MS93}). From \eqref{eq:NUEinf0} and
\eqref{eq:NUEinf1} we obtain
\begin{align}\label{eq:NUEsup}
  \limsup_{n\to+\infty} \frac 1 n \sum_{j= 0}^{n-1} \log
  \Norm{Dg(g^j(x))^{-1}}^{-1} > 0 \quad\text{for Lebesgue
    almost every}\quad x.
\end{align}
Hence, according to \cite{Pinheiro05}, $g$ is a
non-uniformly expanding map since, besides
\eqref{eq:NUEsup}, the orbit of the critical set clearly
\emph{does not accumulate} the critical set.  This ensures
that $g$ admits an absolutely continuous invariant
probability measure, which is unique because $g$ is a
transitive map (see the proof of Lemma~\ref{le:alamane}
below).

To complete the proof of Theorem~\ref{thm:seg}, we are left
to prove Lemma~\ref{le:alamane}.

\begin{proof}[Proof of Lemma~\ref{le:alamane}]
  It is not difficult to see that $g_0$ is topologically
  conjugated to the Tent Map $T(x):=1-2|x|$ for $x\in I$
  under some homeomorphism $h$ of the interval $I$.
  Indeed, searching for $h$ of the form $h(x)=x+u(x)$ on each
  interval $[-1,0]$ and $[0,1]$ for some continuous
  $u:\pm[0,1]\to I$  with small $C^0$-norm, we get the relation
  \begin{align*}
    h(g_0(x))=1-2|h(x)|
    \quad\text{or}\quad
    g_0(x) + u(g_0(x))
    = 1-2|x+u(x)|.
  \end{align*}
  For $x\in [0,1]$ we have $h(x)\ge0$ and so we obtain
  \begin{align*}
    \underbrace{2^{-1}
    u(g_0(x))+u(x)}_{\cL(u)x}=
    \frac12\big(1-2x-g_0(x)\big).
  \end{align*}
  We remark that from this relation it follows that
  $u(0)=u(1)=0$ and, moreover, the right hand side is
  strictly smaller than $1/2$ uniformly on $[0,1]$, that is,
  $\cL(u)\le1/2-\xi$ for some $0<\xi<1/2$.  Clearly
  $\cL(u)=\cI+L$ where $\cI$ is the identity on
  $C^0([0,1],I)$ and $(Lu)x=2^{-1}u( g_0(x) )$ has $C^0$
  norm $\le1/2$. Thus the linear operator $\cL:
  C^0([0,1],I)\to C^0([0,1],I)$ admits an
  inverse. Analogously for the conjugation equation on
  $[-1,0]$. So we can find $h=\cI+u$ with $u$ having $C^0$
  norm $<1$, as needed to ensure that $h$ is invertible,
  thus a homeomorphism of $I$.

  \emph{This guarantees that $g_0$ is transitive and, in
    particular, has no attracting periodic orbits.}

  The subset $K:=\cap_{n\ge0} T^{-n}([-1,1-\epsilon])$ is a
  $T$-invariant Cantor set with zero Lebesgue measure, for
  each $0<\epsilon<1$. Thus $K_0:=h(K)$ is also a
  $g_0$-invariant Cantor set such that, for some $\delta>0$
  and every $z\in K_0$ satisfies $g_0^n(z)\in
  I\setminus\big( (-\delta,\delta)\cup (1-\delta,1] \big)$
  for all $n\ge0$. 

  That is, $K_0$ is the set of points whose future orbits
  under $g_0$ do not visit a neighborhood $V_0$ of the
  critical set, so that $g_0\mid (I\setminus
  \overline{V_0})$ is a $C^\infty$ map acting on $K_0$, by
  the definition of $g_0$.  Hence in $K_0$ we have no
  critical points and no attracting periodic orbits, thus
  the restriction $g_0\mid K_0: K_0\to K_0$ is a uniformly
  expanding local diffeomorphism by Ma\~n\'e's results
  in~\cite{Man85}.

  Therefore the Lebesgue measure of $K_0$ is zero, for
  otherwise this set would have a Lebesgue density point
  $p$. Since $g_0$ is $C^\infty$ is a neighborhood of $K_0$,
  this would imply that $K_0$ would contain some interval
  $J$ (see e.g. \cite{Vi97b} or
  \cite{alves-luzzatto-pinheiro2005}). But $g_0$ is
  uniformly expanding on $K_0$, thus the length of the
  successive images $g_0^k(J)$ of $J$ would grow to at least
  the length of one of the domains of monotonicity of $g_0$,
  in a finite number of iterates. The next iterate would
  contain the critical point, contradicting the definition
  of $K_0$.

  It follows that the set $E(\delta)$ of points of $I$ which
  do not visit a $\delta$-neighborhood of $c_0$ under the
  action of $g_0$ has zero Lebesgue measure, for all small
  $\delta>0$. This ensures that $\cup_{n>N}\cap_{k>n}
  g_0^{-k}\big( E(k^{-1})\big)$ has zero volume for every
  big $N>1$. Consequently the set $\cap_{n>N}\cup_{k>n}
  g_0^{-k}\big(M\setminus E(k^{-1})\big)$ has full measure,
  and for points in this set there are infinitely many
  iterates visiting any given neighborhood of $c_0$.
\end{proof}

This completes the proof of Theorem~\ref{thm:seg}.

\subsection{The unperturbed basic dynamics}
\label{sec:unpert-basic-dynamic}

We now adapt the example in Section \ref{sec:an-example-nue}
in order to obtain a map $f$ which will be a kind of
Poincar\'e return map for a singular flow, that we will
perturb later to obtain a Rovella-like flow.  We again
construct a map $f$ in the torusphere by defining its action
in the meridians and parallels.

Let $f_0: I \to I$ be a $C^{1+}$ non-flat unimodal map with the
critical point $c=0$, that is,
$$
f_0(x)=\begin{cases}
\psi^+(x^\alpha) & \text{if   } x\in(0,1]\\
\psi^-(|x|^\alpha) & \text{if   } x\in[-1,0)
\end{cases};
$$
for smooth monotonous increasing diffeomorphisms
$\psi^\pm:[0,1]\to I$.  Moreover we assume that the critical
order $\alpha$ is at least $2$ and that $f$ satisfies (see
Figure~\ref{fig:1dcritical}):
\begin{figure}[htbp]
\begin{center}
\psfrag{s}{$-1=p_0$}\psfrag{p}{$p_1$}\psfrag{q}{$p_2$}\psfrag{n}{$+1$}
\psfrag{0}{$0$}
  \includegraphics[scale=0.3 ]{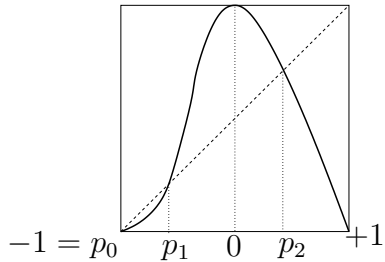}
\end{center}
\caption{The one-dimensional map $f_0$.}
  \label{fig:1dcritical}
\end{figure}
\begin{enumerate}
\item $f_0(\pm1)=-1$, $f_0$ has exactly three fixed points:
  $p_0=-1< p_1 <c=0 < p_2<1$;
\item $f_0'(p_2)<-1<0 \leq f_0'(-1) < 1< f_0'(p_1)$.
\end{enumerate}

\begin{figure}[htpb]
  \centering
  \psfrag{s}{$-1$}\psfrag{n}{$+1$}\psfrag{p}{$\TT^k_1$}
  \psfrag{N}{$N$}\psfrag{u}{$\vec u$}\psfrag{v}{$\vec v$}
  \includegraphics[scale=0.4]{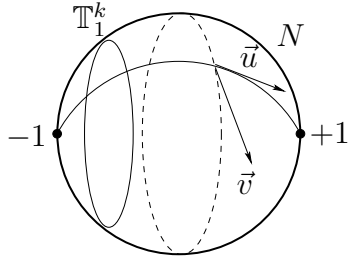}
  \caption{\label{fig:repel}The parallel corresponding to
    the expanding fixed point is a uniform repeller for $f$,
    a vector $\vec u$ tangent to a meridian, and a vector
    $\vec v$ tangent to a parallel.}
\end{figure}

The map $f_0\times
f_1$ induces a $C^{1+\alpha}$ map $f : N\to N$ by
$f=\Upsilon \circ (f_0\times f_1 ) \circ \Upsilon^{-1}.$
Let $\TT_1^k=\Upsilon(\{p_1\}\times \TT^k)$.  Note that $f(\TT_1^k)=
\TT^k_1$, in other words, $\TT^k_1$ is positively invariant
by $f$ and a uniform repeller, see Figure~\ref{fig:repel}.

% Let $\hat\Sigma=\RR^{3k}$. We define
% $F:\cT^k\subset\hat\Sigma\to \hat\Sigma$ by its action on
% each solid torus
% \begin{align*}
%   \cT=B_{1/2}(\sS^1\times 0)=\{ Y\in\RR^3:
%   \dist(Y,\sS^1\times 0)\le 1/2\},
% \end{align*}
% where $\sS^1\times 0=\{(x_1,x_2,0)\in\RR^3:
% x_1^2+x_2^2=1\}$.  On each such torus the map is the Smale
% solenoid map sending the solid torus inside itself (for the
% description of this map see e.g.~\cite{Sm67}). The torus
% $\TT^k$ is the quotient of $\cF^k$ by the stable leaves.

%%%%%%%%%%%%%%%%%%%%%%%%%%%%%%%%%%%%%%%%%%%%%%%%%%%%%%

\subsection{The unperturbed singular flow}
\label{sec:unpert-singul-flow}

Here we build a geometric model for a $(k+5)$-dimensional flow
$X^t_0$, $t \ge 0$.  We write $B^n$ for the $n$-dimensional unit
ball in $\RR^n$, that is  $B^n:=\{ x=(x_1,\dots,x_n)\in\RR^n:
\sum_{i=1}^n x_i^2<1\}$.

Recall that $f_1: \TT^k \to \TT^k$ is the expanding map defined
in section \ref{sec:an-example-nue}. Such map has a lift to an inverse limit
which is a higher dimensional version of a Smale solenoid map,
see \cite{Sm67}.

More precisely, as shown in Appendix \ref{sec:isotopy},
given a solid $k$-torus $\cT:=\TT^k\times\DD$, where
$\DD=\{z\in\CC:|z|\le1\}$, there exists a map $S: e(\cT)\to
e(\cT)$ defined on the image of a smooth embedding $e:\cT\to
B^{k+2}$ such that $\pi_{\DD}\circ S = f_1\circ\pi_\DD,$
where $\pi_\DD:e(\cT)\to e(\TT^k\times 0) \simeq\TT^k$ is
the projection along the leaves of the foliation $\cF^s:=\{
e(\Theta\times\DD) \}_{\Theta\in\TT^k}$ of $\cT$. Moreover,
$F$ contracts the disks in $\cF^s$ by a uniform contraction
rate $\lambda\in(0,1)$.  We remark that $\TT^k$ is the
quotient of $e(\cT)$ over $\cF^s$. From now on we identify
$e(\cT)$ with $\cT$ and $\TT^k$ with $e(\TT^k\times0)$.  The
map $S$ is the higher dimensional version of a Smale
solenoid map we mentioned above.

We also have (see Appendix~\ref{sec:isotopy}) that there
exists an smooth isotopy $\phi_t:B^{k+2}\to B^{k+2}$ between
$S=\phi_1$ and the identity map $\phi_0$ on $B^{k+2}$.

We define a flow between two cross sections
$\Sigma_j=\pi_1(I^2\times\{3-j\}\times B^{k+2}), j=1,2$,
by
$$
Y^t\big( \pi_1 (x_1,x_2,2,W)\big)=
\pi_1\big(x_1,x_2,2-t,\phi_t(W)\big),\quad 0\le t\le1
$$
where $(x_1,x_2,2-t,W)\in
I^2\times[1,2]\times B^{k+2}$ and
$$
\pi_1: \RR\times\RR\times\RR\times B^{k+2}\circlearrowleft,
\quad (x_1,x_2,x_3,W)\mapsto(x_1,x_2,x_3,(1-x_1^2)^{1/2}W).
$$
We note that $\Sigma_j\simeq I\times\overline{B^{k+3}}$
naturally for any $j=1,2$.
We extend this flow so that there exists a Poincar\'e return
map from $\Sigma_2$ to $\Sigma_1$ with the properties we
need. For this we take a linear flow on
$\RR^3\times \RR^{k+2}$ with a singularity $s_0$ at the
origin having real eigenvalues $\lambda_1 >0$, and
$\lambda_1+\lambda_j< 0$ for $2\le j \le k+5$. For
simplicity, we also assume that $\lambda_j=\lambda_3$ for
all $j\ge4$ and that $\alpha:=-\lambda_3/\lambda_1$ and
$\beta:=-\lambda_2/\lambda_1$ satisfy $\beta>\alpha+2$. This
last \emph{strong dissipative condition} on the saddle $s_0$
ensures that the foliation corresponding to the $x_2$
direction is dominated, and so persists for all $C^2$ nearby
flows.

We note that the subspace $\{0\}\times\RR^2\times B^{k+2}$,
excluded from further considerations, is contained in the stable
manifold of $s_0$ and so its points never return no
$\Sigma_j$.

\begin{figure}[htbp]
\begin{center}
\psfrag{S2}{$\Sigma_1$}\psfrag{S1}{$\Sigma_2$}
\psfrag{Y}{$Y^1$}\psfrag{p}{$s_0$}
\psfrag{P1}{$\pi_1$}\psfrag{P0}{$\pi_2$}\psfrag{GP2}{$P_2$}
\psfrag{GP1}{$P_1$}\psfrag{GP0}{$P_0$} \psfrag{N}{$N$}
\psfrag{T}{$T^+$}\psfrag{t}{$T^-$}\psfrag{GP3}{$P_3$}
  \includegraphics[scale=0.4]{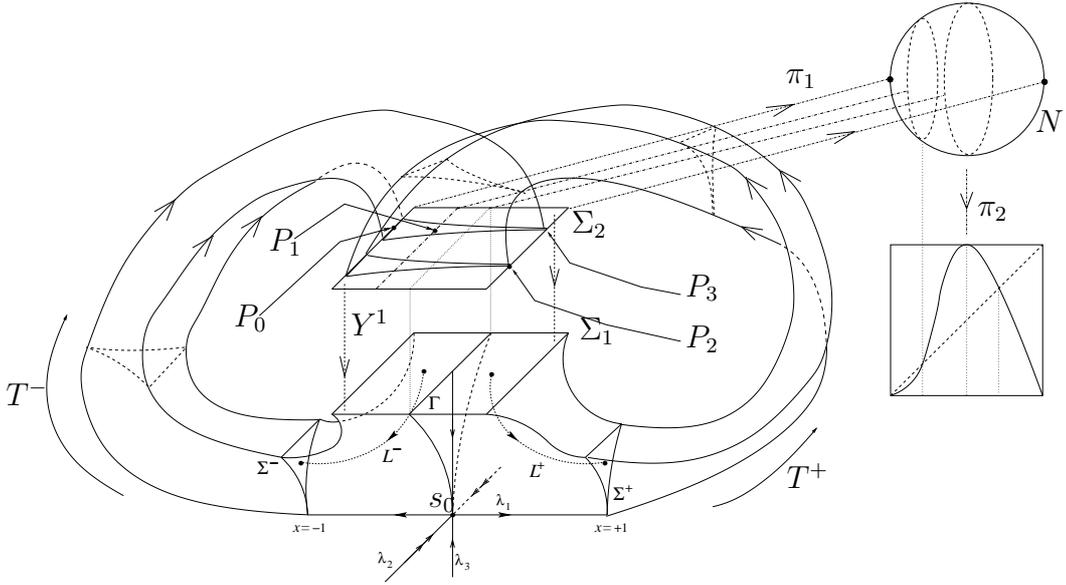}
\end{center}
\caption{The starting singular flow with the transformations
  $Y^1$, $L^\pm$ and the projections $\pi_1$ and $\pi_2$.}
  \label{fig:startflow}
\end{figure}

We write $(x_1,x_2,1,W)$ for a point on $\Sigma_2$, with
$x_i\in\RR$ and $W\in B^{k+2}$, and consider the
cross-sections $\Sigma^\pm=\{(\pm1,x_2,x_3,W): (x_2,x_3)\in
I^2, W\in\RR^{k+2}\}$ to the flow and the Poincar\'e first
entry transformations given by (see
Figure~\ref{fig:startflow}):
\begin{align}\label{eq:thrusing}
  L^\pm: \Sigma_1\cap\{\mp x_1>0\}\to \Sigma^\pm: \quad
  (x_1,x_2,1,W)\mapsto (\pm 1, x_2 |x_1|^\beta,
  |x_1|^{\alpha}, |x_1|^{\alpha}W).
\end{align}

\begin{remark}
  \label{rmk:Holder-thru}
  The maps $L^\pm$ given in \eqref{eq:thrusing} are clearly
  H\"older maps in their domain of definition. Moreover the
  time the flow needs to take the point $(x_1,x_2,1,W)$ to
  $\Sigma^\pm$ is given by $-\log |x_1|$, where $|x_1|$ is
  the distance to the local stable manifold of $s_0$ on
  $\Sigma_2$.
\end{remark}

Now we define diffeomorphisms from the image of $L^\pm$ to
$\Sigma_1$ which can be realized as the first entry maps
from $\Sigma^\pm$ to $\Sigma_1$ under a flow defined away
from the origin in $\RR^{k+5}$. Since we want to define a
flow with an attractor containing the origin in $\RR^{k+5}$, we
need to ensure that the return map defined on $\Sigma_1$
through the composition of all the above transformations
does preserve the family of tori together with their stable
foliations. Moreover the quotient of this return
transformation over the stable directions should be the map
$f$.  We write these transformations as
$T^\pm:\Sigma^\pm\to\Sigma_2$ given by
$$
(\pm1,z_2,z_3,V)
\mapsto\left(\psi^\pm(z_3),
 \pm\frac12+\frac{z_2}C,
2,
\Psi^\pm(z_3)V\right)
$$
where $C>0$ is big enough so that the map restricted to the
first $3$ coordinates is injective and thus a diffeomorphism
with its image. We recall that $\psi^\pm$ is part of the
definition of the one-dimensional map $f_0$.  The
diffeomorphisms $\Psi^\pm:[0,1]\to I$ are chosen to ensure
that the quotient map is well defined on $N$, as follows:
$$
z_3\mapsto
\Psi^\pm(z_3):=
\frac1{z_3}\sqrt{\frac{1-\psi^\pm(z_3)^2}{1-z_3^{2/\alpha}}},
\quad\text{for } z_3\in(0,1)
$$
and we set $\Psi^\pm(0):=\mp1$ and $\Psi^\pm(1)=\pm1$. We
remark that $\Psi^\pm$ are diffeomorphisms, see
Figure~\ref{fig:startflow}.  Indeed both the numerator and
denominator inside the square root can be written as a
Taylor series around $z_3=1$ with an expression like
$
const\cdot(1-z_3) + o(1-z_3)
$
for some non-zero constant, thus $\Psi^\pm$ is
differentiable at $1$.  Now $\psi^\pm$ can be expanded
around $z_3=0$ as
$
1+const\cdot z_3 + o(|z_3|)
$
for some positive constant. Hence $\sqrt{1-(\psi^\pm)^2}$
can be expanded as
$
\big(1-(1-const\cdot z_3 + o(|z_3|))^2\big)^{1/2}=
\const\cdot z_3 + o(|z_3|)
$
thus $\Psi^\pm$ is also differentiable at $0$.

Now we check that the return map $R_0$ given by
$
R_0:= T^\pm \circ L^\pm \circ Y^1:\Sigma_2\circlearrowleft
$
can be seen as a map on $N$.  Indeed notice that $Y^1$
commutes with $\pi_1$ by construction and that for $x_1>0$
\begin{align*}
  (T^+\circ L^+\circ \pi_1)(x_1,x_2,1,W) &= (T^+\circ
  L^+)\big(x_1,x_2,1,(1-x_1^2)^{1/2} W\big)
  \\
  &= T^+\Big( \psi^+\big(|x_1|^\alpha\big), x_2
  |x_1|^\beta, |x_1|^\alpha, |x_1|^\alpha
  (1-x_1^2)^{1/2} W\Big)
  \\
  &= \left( f_0(x_1), \frac12+\frac{x_2
      |x_1|^\beta}C , 2, W
    \sqrt{1-\psi^+(|x_1|^\alpha)^2} \right)
  \\
  &= \pi_1\Big( f_0(x_1), \frac12+\frac{x_2
    |x_1|^\beta}C , 2, W\Big).
\end{align*}
Hence we  get
\begin{align}
R_0(x_1,x_2,2,W)&= (Y^1\circ T^+\circ L^+ \circ
\pi_1)(x_1,x_2,2,W)
\nonumber
\\
&= \pi_1
\Big( f_0(x_1), \frac12+\frac{x_2
  |x_1|^\beta}C ,
2, \phi_1(W)\Big), \label{eq:returnmap}
\end{align}
where $x_1>0$ and $\phi_1=F$ by the definition of $Y^1$.
Analogous calculations are valid for $x_1<0$ taking the maps
$T^-$ and $L^-$ into account. By the definition of $\phi_t$
we get that $R_0$ maps $\hat \Sigma =
\Sigma_2\cap\pi_1(I^2\times\{1\}\times\cT)$ inside itself.

\begin{remark}
  \label{rmk:domination}
  The contracting direction along the eigendirection of the
  eigenvalue $\lambda_2$ can be made dominated by the other
  directions in this construction by increasing the
  contraction rates given by $C$ and $\lambda_2$. This is
  very important to ensure the persistence of the stable
  lamination, see e.g. \cite{HPS77} and
  Section~\ref{sec:x-chaotic-with}. This is why we assumed
  the strong dissipative condition $\beta>\alpha+2$ on the
  saddle equilibrium $s_0$.  We note that this domination is
  for the action of the flow we are constructing.
\end{remark}

We now remark that quotienting out the contracting
directions of $\cF^s$ we obtain the map $f$ on $N$.  Indeed,
each element of the quotient can be seen as the image of the
following projection
$
\pi_2:\hat\Sigma_2\to N$ given by
$(x_1,x_2,1,W)\mapsto(x_1,(1-x_1^2)^{1/2}\pi_\DD(W))$
where, we recall, $\pi_\DD: \cT\to\TT^k=\sS^1 \times
\overset{k}{\dots} \times \sS^1 \subset \CC^k \simeq \RR^{2k}$.
% $$
% \big((x_1,y_1,z_1),\dots,(x_k,y_k,z_k)\big)\mapsto \left(
%   \frac{(x_1,y_1)}{\sqrt{x_1^2+y_1^2}}, \dots,
%   \frac{(x_k,y_k)}{\sqrt{x_k^2+y_k^2}} \right).
% $$
It is easy to see that \emph{the return map $R_0$ is
  semiconjugated to $f$ through $\pi_2$}, that is $
\pi_2\circ R_0 = f\circ \pi_2.  $

\begin{remark}
  \label{rmk:cylinder}
  We take advantage of the fact that the hyperplanes
  $\{x_1=\pm1\}$ can be identified to a single point due to
  the dynamics of $\hat g_0$. If we do not perform this type
  of identification, i.e. if we consider instead the
  projection $\pi_3:\hat\Sigma_2\to I\times\TT^k$ given by
  $(x_1,x_2,1,W)\mapsto(x_1,\pi_\DD(W))$,
  then we get a cylinder $I\times\TT^k$ as the domain of the
  quotient map, instead of the torusphere $N$, and likewise
  $\pi_3\circ R_0 = (f_0\times f_1)\circ \pi_3$, where
  $f_1:\TT^k\circlearrowleft$ is the expanding map on the
  $k$-torus defined in Section~\ref{sec:an-example-nue}.
\end{remark}

\subsubsection{Localizing some periodic orbits of the flow}
\label{sec:some-period-orbits}

The vector field $X_0$ just constructed has a flow with an
attracting periodic orbit, the orbit of
$P_0=\pi_1(-1,y_*,2,0)$ for some $y_*\in I$. Indeed, we note
that $-1$ is an attracting fixed point for the
one-dimensional map $f_0$ and that the $f_0$ orbit of almost
every $x\in I$ tends to $-1$. Then, since on the second
coordinate in $\Sigma_2$ we have a strong contraction under
the action of $R_0$, this ensures that there exists $y_*$ as
above satisfying $R_0(\pi_1(-1,y_*,2,0))=\pi_1(-1,y_*,2,0)$.
We note that along the $x_1$ and $x_2$ directions the flow
clearly is a contraction.  Moreover, for $x$ in the interval
$(-1,p_1)$, the "toruspherical coordinates" of
$\pi_1(x,y,2,W)$ (that is, the last coordinate of dimension
$k+2$) tend to $0$ since they are multiplied by
$\sqrt{1-f_0^k(x)^2}$ and $f_0^k(x)\to-1$ as $k\to+\infty$.

In addition the fixed point $p_1$ of $f_0$ corresponds to a
hyperbolic invariant subset for the flow inside the solid
torus $P_1=\pi_1((p_1,y^*,2)\times\cT)$ for some $y^*\in I$.
We observe that quotienting out the stable directions we get the
invariant torus $\TT^k_1$ as already mentioned before, see
Figure~\ref{fig:startflow}.

Finally the orbit of $P_2=\pi_1(1,\bar y,2,0)$ is in the
stable manifold of $P_0$ for any $\bar y\in I$. Indeed it
returns to the stable leaf $\{\pi_1(-1,y,2,0):y\in I\}$ and
never leaves this leaf in all future returns. In the
quotient $\pi_2(\hat\Sigma_2)$ the point $P_0$ is fixed and
$P_2$ is one of its preimages under $f$.

%%%%%%%%%%%%%%%%%%%%%%%%%%%%%%%%%%%%%%%%%%%%%%%%%%%%%%%

\subsection{Perturbing the original singular flow}
\label{sec:pert_flow}

Now we make a  perturbation $X$ of the vector field $X_0$
constructed in the previous section to obtain a Rovella-like
attractor.

From now on we assume that $P_2$ is the point where the
component of the unstable manifold of $s_0$ through
$\Sigma^+$ first arrives at $I^2\times\{2\}\times B^{k+2}$.
We consider then the positive orbit of $P_2$ up until it returns to
$\hat\Sigma_2$.  By construction the return point $R_0(P_2)$
has the expression $\hat P=\pi_1(-1,\hat y,2,0)$, so it
returns in $\pi_2(-1,0)$, after the quotient through
$\pi_2$.

\begin{figure}[htpb]
  \centering
  \psfrag{0}{$0$}\psfrag{p}{$p_1$}\psfrag{r}{$p_0=-1$}
\psfrag{s1}{$s_1$}\psfrag{P0}{$\pi_2$}
\psfrag{N}{$N$}\psfrag{n}{$+1$}\psfrag{s}{$-1$}
\psfrag{T}{$\TT_1^k$}\psfrag{sr}{$\hat s$}
  \includegraphics[scale=0.5]{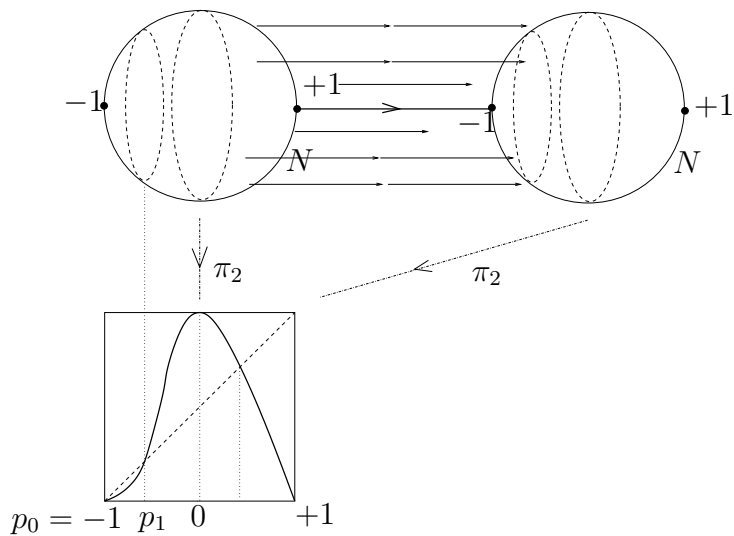}
  \caption{\label{fig:beforeperturb} This represents the
    initial flow before the perturbation and the
    one-dimensional quotient.}
\end{figure}

This is a regular orbit, so it admits a tubular
neighborhood. We can assume that the orbit chosen above
returns to $\hat\Sigma_2$ close enough to $P_0$, that is $\hat
y$ is close to $y_*$. We also assume that $P_1$ is also very
close to $P_0$ so that the tubular neighborhood contains
both $P_0$ and $P_1$, see
Figure~\ref{fig:beforeperturb}.

We note that since $P_0$ is a hyperbolic attracting orbit it
is easy to extend the tubular neighborhood to its local
basin of attraction whose topological closure contains
$P_1$.

\begin{figure}[htbp]
\begin{center}
\psfrag{0}{$0$}\psfrag{p}{$p_1$}\psfrag{r}{$p_0=-1$}
\psfrag{s1}{$s_1$}\psfrag{P0}{$\pi_2$}\psfrag{sr}{$\hat s$}
\psfrag{N}{$N$}\psfrag{n}{$+1$}\psfrag{s}{$-1$}
\psfrag{T}{$\TT_1^k$}
  \includegraphics[scale=0.6]{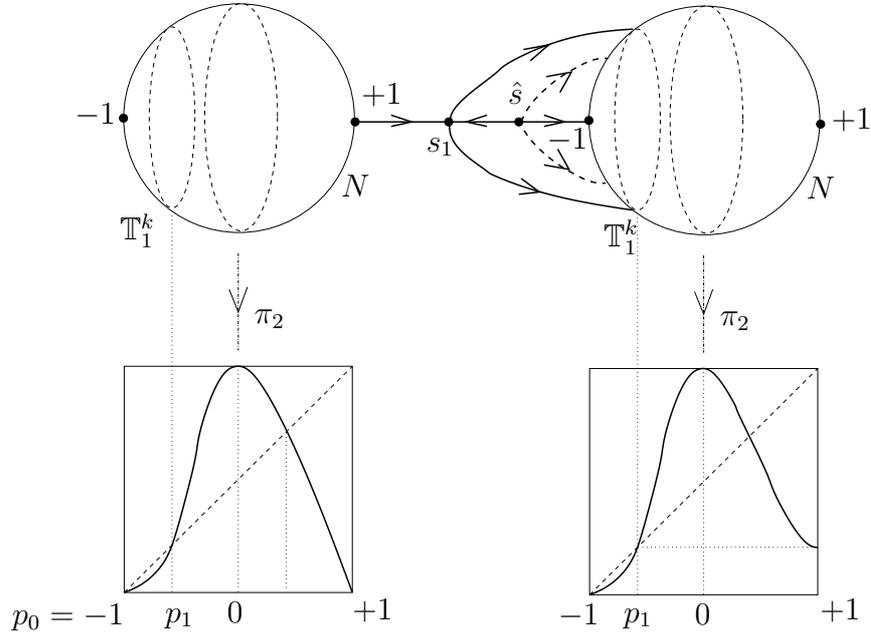}
\end{center}
\caption{The perturbed singular flow seen through the
  projections $\pi_1$ and $\pi_2$.}
  \label{fig:tubularperturb}
\end{figure}

In this setting we can now perturb the flow inside the
tubular neighborhood in the same way as to produce a kind of
Cherry flow, introducing two hyperbolic saddle singularities
$\hat s$ and $s_1$. Here the extra dimensions are very
useful to enable us to introduce such saddle fixed point
with the adequate dimensions of stable and unstable
manifolds.  The saddle $s_1$ has $k+1$ expanding eigenvalues
$\bar\lambda_0,\dots,\bar\lambda_{k}$, the remaining $4$
contracting eigenvalues $\bar\lambda_j, j=k+1,\dots,k+4$,
and $s_1$ is sectionally dissipative:
$\bar\lambda_i+\bar\lambda_j<0$ for all $i\le k$ and $j>k$.
In fact the extra dimensions are essential to allow the
construction of such a saddle.

The other saddle $\hat s$ has $k+2$ expanding eigenvalues
and $3$ contracting ones.  We assume that this perturbation
is done in such a way that the $(k+1)$-dimensional unstable
manifold of $s_1$ \emph{contains} $\TT^k_1$, see
Figure~\ref{fig:tubularperturb}, and the stable manifold of
$s_1$ is everywhere tangent to the subspace given by the
direction of the stable manifolds of the solid tori together
with the $x_2$ direction. In this way we can still quotient
out the stable leaves, which are preserved by the perturbed
flow. \emph{This is the flow of $X$ in the statement of
Theorem~\ref{mthm:pri}.}

\begin{figure}[htbp]
\begin{center}
\psfrag{W}{$W^s(s_1)$}\psfrag{w}{$W^u(s_0)$}
\psfrag{r}{$\hat  s$}\psfrag{V}{$W^{ss}(s_0)$}
\psfrag{s}{$s_1$}\psfrag{t}{$s_0$}
\psfrag{u}{$W^u(s_1)$}\psfrag{v}{$W^s(s_0)$}
  \includegraphics[scale=0.8]{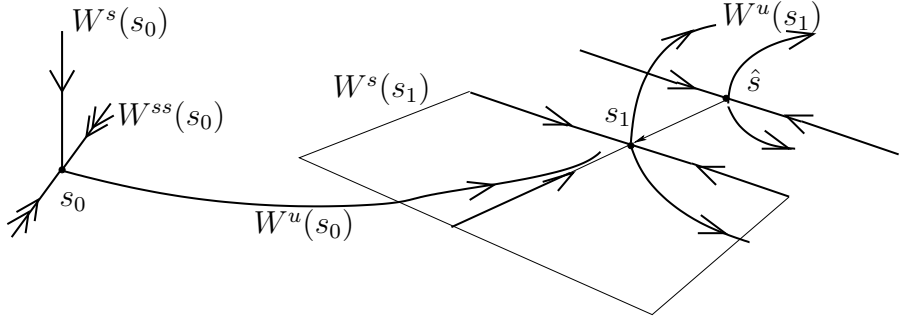}
\end{center}
\caption{The connection between $s_0$ and $s_1$.}
  \label{fig:saddle-connection}
\end{figure}

\begin{remark}\label{rmk:keepsaddleconnection}
  Since $P_2$ is part of the unstable manifold of $s_0$,
  then we have constructed a saddle connection between $s_0$
  and $s_1$. Because the stable manifold of $s_1$ is
  $3$-dimensional, we can keep the connection for nearby
  vector fields restricted to a $(k+2)$-codimension 
  submanifold of the space of all smooth vector fields: all
  we have to do is to keep one branch of the one-dimensional
  unstable manifold of $s_0$ contained in the
  $3$-dimensional local stable manifold of $s_1$, and this
  submanifold has codimension $k+2$. See
  Figure~\ref{fig:saddle-connection}.
\end{remark}

The action of the first return map $R$ of the new flow $X^t$
to $\Sigma_2$, on the family of stable leaves which project
to the interval $[p_1,1]$, equals (except for a linear
change of coordinates) the map $g$ presented in
Section~\ref{sec:an-example-nue}, see
Figure~\ref{fig:tubularperturb}. That is, we have a
multidimensional non-uniformly expanding transformation as
the quotient of the first return map.  This transformation
sends each $k$-torus in $N$ to another $k$-torus and so it
can be further reduced to the interval map $g_0$, by
considering its action on the tori.

\begin{remark}
  \label{rmk:poincare-time}
  The introduction of the saddle-connection changes the
  return time of the points in $\Sigma_2$ under the new flow,
  but we can assume that the distance from the stable
  manifold of $s_1$ of the orbit of
  $(x_1,x_2,2,W)\in\Sigma_2$ near $s_1$ is given by the
  $x_3$ coordinate of $L^\pm(x_1,x_2,2,W)$, that is, by
  $|x_1|^\alpha$. Therefore the Poincar\'e return time is
  now bounded, after Remark~\ref{rmk:Holder-thru}, by some
  uniform constant plus $const\cdot\log|x_1|$. The value
  $|x_1|$ equals the distance of $(x_1,x_2,2,W)$ to the
  critical set of $g$.
\end{remark}

This reduction to a one-dimensional model will be essential
to our analysis of the existence of a absolutely continuous
invariant probability measure for $g$ and the existence of
a physical probability measure for the attracting set of the
flow near the origin.

\subsubsection{The attracting set and an invariant stable foliation}
\label{sec:attract-set-its}

We note that the set $\widetilde\Lambda_\Sigma:=\cap_{n>0}
R^n(\Sigma_2)$ contains a subset $\Lambda_\Sigma:=\cap_{n>0}
R^n(\widetilde U)$, where $\widetilde U$ is a small
neighborhood in $\Sigma_2$ of $\pi_2^{-1}[p_1,1]$, which is
an attracting set for $R$, that is, $\overline{R(\widetilde
  U)}\subset\widetilde U$. 
Hence the saturation
$\Lambda=\cup_{t\in\RR}X^t(\Lambda_\Sigma)$ is
also an attracting set for $X^t$ such that
$\Lambda_\Sigma= \Lambda\cap\Sigma_2$.

It is easy to see that every point $z=\pi_1(x_1,x_2,1,W)$ of
$\Lambda_\Sigma$ belongs to the solenoid attractor in
$\{(x_1,x_2,1)\}\times B^{k+2}$. Therefore it is
straightforward to define a local stable foliation
$\cF^{ss}$ through the points of $\Lambda_\Sigma$: we define
$\cF^{ss}_z$ to be the stable disk through $z$ of the solid
torus $\{(x_1,x_2,1)\}\times \phi_1(\cT)$.  Hence we get a
$DR$-invariant continuous and uniformly contracting
subbundle $E^{ss}$ of $T\Lambda_\Sigma$ given by
$E^{ss}_z=T_z\cF^{ss}_z$.

From the existence of the invariant contracting subbundle
$E^{ss}$ over $T\Lambda_\Sigma$, we can define the normal
subbundle
$G=(E^{ss})^\perp\cap \RR\times\{(0,0)\}\times\RR^{k+2}$
to $E^{ss}$ in the tangent space to $I\times\{(x_2,1)\}\times B^{k+2}$,
and use this pair of continuous bundles to define a stable
cone field and the complementary cone field
\begin{align*}
  C^s(z)=\{(u,v)\in E^{ss}_z\oplus G_z : \|u\|\ge\|v\|\}
,\quad
  C^u(z)=\{(u,v)\in E^{ss}_z\oplus G_z : \|u\|\le\|v\|\}
\end{align*}
for all points $z\in\Lambda_\Sigma$. 

Is is clear that the bundle of tangent spaces to $\{I\times\{(x_2,1)\}\times
B^{k+2}\}_{x_2\in I}$ is invariant under $DR$, under both forward and
backward iteration. Moreover, we have that $C^s$ is strictly
invariant under backward iteration by $DR$ and,
consequently, the complementary cone is invariant under
forward iteration by $DR$. Therefore vectors in $C^u(z)$ make
an angle with vectors in $E^{ss}_z$ uniformly bounded away
from zero. We use this in the next section to prove the
existence of a dominated splitting over $T\Lambda_\Sigma$.

%%%%%%%%%%%%%%%%%%%%%%%%%%%%%%%%%%%%%%%%%%%%%%%%%%%

\section{Properties of the vector field and its unfolding}
\label{sec:conseq-unfold-x}

We now prove all items of Theorem~\ref{mthm:pri} and show
that the flow of $X$ is chaotic with multidimensional
expansion. We also consider its perturbations and prove
Theorem~\ref{mthm:codimension}. Some technical points are
postponed to Sections~\ref{sec:higher-dimens-misiur}
and~\ref{sec:full-basin-attract} and
Appendix~\ref{sec:hyperb-neighb-constr}.

\subsection{$X$ is chaotic with multidimensional nonuniform
  expansion}
\label{sec:x-chaotic-with}

We observe that, since we have a one-dimensional quotient
map where the critical point is mapped to a repelling fixed
point, we are in the setting of ``Misiurewicz maps'', see
the right hand side of Figure~\ref{fig:tubularperturb}.

\begin{figure}[htpb]
  \centering
  \psfrag{0}{$0$}\psfrag{p}{$p_1$}\psfrag{r}{$p_0=-1$}
\psfrag{s1}{$s_1$}\psfrag{P0}{$\pi_2$}
\psfrag{N}{$N$}\psfrag{n}{$+1$}\psfrag{s}{$-1$}
\psfrag{T}{$\TT_1^k$}
  \includegraphics[scale=0.5]{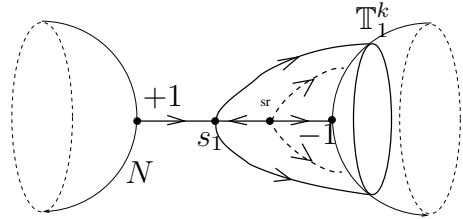}
  \caption{\label{fig:symperturb0} The last perturbation
    defining the flow $X$. We can project the dynamics into
    a one-dimensional map.}
\end{figure}

We note that the map on the torusphere for the perturbed
flow is \emph{non-uniformly expanding in all directions},
because the non-uniformity is seen on every direction away
from the repelling torus $\TT_1^k$, and because the
singularity contracts also in every direction, see
Figure~\ref{fig:symperturb0}. This was already proved in
Section~\ref{sec:an-example-nue}.

\subsubsection{Invariant probability measure}
\label{sec:invari-probab-measur}

The results in Section~\ref{sec:higher-dimens-misiur} ensure
that there exists an ergodic absolutely continuous invariant
probability measure $\upsilon$ for the transformation $g$,
which is the action of the first return map $R$ of the flow
of $X$ to $\Sigma_2$ on the stable leaves. Indeed $g$
satisfies all the conditions of Theorem~\ref{thm:Misiurew}
since its quotient $g_0$ is a Misiurewicz map of the
interval.  Moreover $\upsilon$ is an \emph{expanding
  measure}: all Lyapunov exponents for $\upsilon$-a.e. point
are strictly positive, as shown in
Section~\ref{sec:an-example-nue}.

From standard arguments using the uniformly contracting
foliation through $\Lambda_\Sigma$ (recall
Section~\ref{sec:attract-set-its}), see e.g. \cite[Section
6]{APPV}, we can construct an $R$-invariant ergodic
probability measure $\nu$, whose basin has positive Lebesgue
measure on $\Sigma_2$, and which projects to $\upsilon$
along stable leaves. Finally, using a suspension flow
construction over the transformation $R$ we can easily
obtain (see \cite[Section 6]{APPV} for details) a
corresponding ergodic physical probability measure $\mu$ for
the flow of $X$, which induces $\nu$ as the associated
$R$-invariant probability on $\Sigma_2$.

We claim that $\mu$ is a hyperbolic measure for $X^t$ with
$k+1$ positive Lyapunov exponents. To prove this, we first
obtain a dominated splitting for the tangent $DR$ of the
return map $R$ which identifies the bundle of directions
with nonuniform expansion.

%%%%%%%%%%%%%%%%%%%%%%%%%%%%%%%%%%%%%%%%%%%%%%%%%%%%%%%

\subsubsection{Strong domination and smooth stable foliation}
\label{sec:dominat-splitt}

We now consider the return map $R$ to $\Sigma_2$ and show
that the subbundle $E^{s}$ corresponding to the tangent
planes to the local stable leaves of the solenoid maps is
``strongly $\gamma$-dominated'' by the non-uniformly
expanding direction ``parallel'' to the toruspheres, for
some $\gamma>1$.

Let $\cF^{s}_{a,b}$ be the bidimensional stable foliation
of the solenoid map acting on the section
$\{x_1=a,x_2=b\}\cap (I^2\{1\}\times B^{k+2})$ for
each $(a,b)\in I^2$. Then, from the expression
\eqref{eq:returnmap} of $R_0$, which is only modified to $R$
by putting $g_0$ in the place of $f_0$, we see that
the derivative of $R_0$ along some leaf
$\gamma\in\cF^{s}_{a,b}$ equals the derivative of $\phi_1$
along the same direction multiplied by $(1-g_0(a)^2)^{1/2}$.
Since $\|D\phi_1\mid\gamma\|\le\lambda$ for some constant
(not dependent on $a,b$ nor on the particular leaf of
$\cF^{s}_{a,b}$) $\lambda\in(0,1)$, we get that
$\|DR\mid\gamma\|\le\lambda(1-g_0(a)^2)^{1/2}.$

We can estimate expansion/contraction rates of the
derivative of the quotient map $g$ as in
Section~\ref{sec:an-example-nue}. Since this map is obtained
through a projection of $R$ along $\cF^{s}$, the real
expansion and contraction rates of the derivative of $R$
along any direction in the complementary cone $C^u$ are
bounded by the corresponding rates of $Dg$ up to
constants. These constants depend on the angle between the
stable leaves and the direction on the complementary cone
field, which is uniformly bounded away from zero.

Hence, recalling that $m(x)$ is the minimum norm  of
$Dg(x)$ for $x=\pi_1(t,b,1,z,\Theta)$ with $t\in I$,
$z\in\DD$ and $\Theta\in[0,2\pi)^k$, to obtain the
smoothness of the foliation, it is enough to get that
\begin{align}
  \label{eq:domination}
  d(x):=\frac{(1-g_0(t)^2)^{1/2}}{\|Dg(x)\|^\gamma m(x)^\omega}
\end{align}
is bounded by some constant uniformly on every point $x$, for some
$\gamma>1$ and $\omega>0$. This implies that for a small
enough $\lambda>0$ we have
$  \lambda(1-g_0(t)^2)^{1/2}\|Dg(x)\|^{-\gamma}<
  \lambda m(x)^\omega < m(x)$,
since $\omega>0$ and $m(x)$ is  bounded from above.
\emph{This ensures (see \cite[Theorem
6.2]{HP70}) that $\cF^{s}$ is a $C^\gamma$ foliation, so
that holonomies  along the leaves of $\cF^{s}$ are of class
$C^\gamma$.}

From Section~\ref{sec:an-example-nue} we know that
\begin{align}
  \label{eq:minimum}
  m(x)&=m(t)=\min\left\{|g_0^\prime(t)|,
    2 \Big|\frac{\cos(\frac\pi2 g_0(t))}{\cos(\frac\pi2 t)}\Big|
  \right\},\quad\text{and}
  \\
  \|Dg(x)\|&=\max\left\{|g_0^\prime(t)|,
    2 \Big|\frac{\cos(\frac\pi2 g_0(t))}{\cos(\frac\pi2 t)}\Big|
  \right\}.
\end{align}
Hence $d(x)=d(t)$ only depends on $t\in I$.  

Now we note that for $t\in I\setminus\{-1,0,1\}$ the
quotient $d(t)$ is continuous. Therefore, if we show that
$d$ can be continuously extended to $\{-1,0,1\}$, then $d$
is bounded on $I$ and $\lambda d(t)$ can be made arbitrarily
small letting the contraction rate $\lambda$ be small
enough, which can be done without affecting the rest of the
construction. Having this concludes the proof of the
smoothness of $\cF^{s}$.

Finally, we compute, on the one hand
$\lim_{t\to\pm1}|{\cos(\frac\pi2
      g_0(t))/\cos(\frac\pi2 t)}|
  =|g_0^\prime(\pm1)| \neq0
$
which shows that $d$ can be continuously extended to $\pm1$.
On the other hand, by the choice of the map $g_0$ in Section
\ref{sec:an-example-nue}, we have that both $m(x)$ and
$\|Dg(x)\|$ are of the order of $|t|^{\alpha-1}$ for $t$ near
$0$ (ignoring multiplicative constants).
Thus
\begin{align}\label{eq:crucial}
  d(x)
  =O\left(\frac{(|t|^\alpha)^{1/2}}{(|t|^{\alpha-1})^\gamma
      (|t|^{\alpha-1})^\omega}\right)
  =O(|t|^{\alpha/2-(\gamma+\omega)(\alpha-1)}).
\end{align}
For $1<\alpha<2$ we have $\alpha/(2(\alpha-1)) > 1$ so that,
in this setting, we can take $\gamma>1$ and $\omega>0$ in
order that $\alpha/2-(\gamma+\omega)(\alpha-1)\ge0$.  Hence
$d$ can also be extended continuously to $0$. This concludes
the proof that $\cF^{s}$ is strongly dominated by the
action of $DR$ along the directions on the complement $C^u$
of the stable cone field, so that it becomes a $C^\gamma$
foliation.

Now we define the subspace $E$ to be the sum of the tangent
space $E^{s}:=T\cF^{s}$ to $\cF^{s}$ with the $x_2$
direction $E^{ss}$ on $\Sigma_2$, i.e. $E:= E^{ss}\oplus E^s$.
This bundle $E$ is a $DR$-invariant contracting subbundle
which is also strongly dominated by the directions on the
center-unstable cone $C^u$ (see
Section~\ref{sec:attract-set-its}), since the direction
$E^{ss}$ can be made even more strongly contracted than the
bundle $E^{s}$, see Remark~\ref{rmk:domination}. In fact,
the same argument as above, especially the
relation~\eqref{eq:crucial}, is analogous. The foliation
$\cF$ tangent to $E$ is then uniformly contracting, with
$3$-dimensional $C^1$-leaves and $C^\gamma$ holonomies, for
some $\gamma>1$.

\begin{remark}
  \label{rmk:subfoliation}
  The foliation $\cF^{ss}$ is a subfoliation of
  $\cF$ in the sense that every leaf
  $\gamma\in\cF$ admits a foliation
  $\gamma\cap\cF^{ss}$ by leaves of $\cF^{ss}$ tangent to
  $E^{ss}$ at every point.
\end{remark}

\begin{remark}
\label{rmk:1dfoliation}
\emph{Arguing with the flow of $X$}, there exists a
one-dimensional foliation $\cF^{ss}$ tangent to the
one-dimensional field of directions
$\cup_{t\ge0}DX^t(E^{ss}_z)$, since this line bundle is
uniformly contracted by $X^t$.  This one-dimensional bundle
is also strongly dominated by the ``saturated'' bundle
$\cup_{t\ge0}DX^t(E^{s}\oplus E^X)$, by the choice of the
constants $C$ and $\beta$ in the construction of the vector
field, see Section~\ref{sec:unpert-singul-flow}.
\end{remark}

We can now complete the construction of the center subbundle
$E^c$. The domination just obtained shows that the
complementary cone field $C^u$ through the points of
$\Lambda_\Sigma$ is strictly invariant by forward iteration
under $DR$, so there exists a unique $DR$-invariant
subbundle $E^c$ contained in $C^u$ and defined on all points
of $\Lambda_\Sigma$. We thus obtain a dominated splitting
$E^s\oplus E^c$ of the tangent bundle of $\Sigma_2$ over
$\Lambda_\Sigma$.

%\margem{da pra estender a decomposicao a todo ponto de $\Lambda$?}

\subsubsection{Hyperbolicity of the physical measure}
\label{sec:hyperb-physic-measur}

The suspension $\mu$ of the ergodic and physical invariant
probability measure $\nu$ for $R$ is also an ergodic and
physical measure for $X^t$ on $U$. In addition, denoting
$\tau(z)$ the Poincar\'e return time for $z\in\Sigma_2$
(which a well defined smooth function except on $\{0\}\times
I\times\{1\}\times B^{k+2}$) and
$\tau^n(z)=\tau(R^{n-1}(z))+\dots +\tau(z)$ for all
$n\in\ZZ^+$ such that $\tau^n(z)<\infty$, we have
$DR^n(z)=P_{R^n(z)}\circ DX^{\tau^{n}(z)}\mid T_z\Sigma_2$,
where $P_z:T_zM\to T_z\Sigma_2$ is the projection parallel
to the direction $X(z)$ of the flow at
$z\in\Sigma_2$. Therefore we can write, for $z\in\Sigma_2$
such that $R^n(z)$ is never in the local stable manifold of
$s_0$ for $n\in\ZZ^+$ and $v\in E^c(z)\subset T_z\Sigma_2$
\begin{align}\label{eq:exponent}
  0<\limsup_{n\to+\infty}\frac1n\log\|DR^n(z)v\| \le
  \limsup_{n\to+\infty}
  \frac{\tau^n(z)}n\cdot\frac1{\tau^n(z)}\log\|DX^{\tau^n(z)}(z)v\|.
\end{align}
But $\tau^n(z)/n\xrightarrow[n\to\infty]{}\nu(\tau)$ is
finite for $z$ in the basin of $\nu$. Indeed, by
Remark~\ref{rmk:poincare-time}, $\tau(z)$ is essentially the
logarithm of the distance to the critical set of $R$,
i.e. the intersection of the local stable manifold of $s_0$
with $\Sigma_2$. Indeed, the main contribution to the return
time comes from the time it takes $z\in\Sigma_2$ to pass
near the singularities, which is given by the logarithm of
the distance to $\Sigma_2\cap W^s_{loc}(s_0)$ and this same
value controls the time it takes the point to pass near
$s_1$ also, except for a multiplication by a positive
constant, because of the local expression of the flow near a
hyperbolic equilibrium and by the form of the connection
between $s_0$ and $s_1$, see
Figure~\ref{fig:saddle-connection} and
Remark~\ref{rmk:poincare-time}.

In the quotient dynamics of $R$, i.e. for the map $f$ on
$N$, the function $\tau$ is comparable to the logarithm of
the distance to the critical set.
This function is $\nu$-integrable as a
consequence of the non-uniformly expanding properties of the
map $f$, as stated in Theorem~\ref{thm:Misiurew} of
Section~\ref{sec:higher-dimens-misiur}.

Hence \eqref{eq:exponent} implies that the Lyapunov
exponents along $\mu$-almost every orbit of $X^t$ are
positive along the directions of the bundle
$E^{cu}(X^t(z)):=DX^t(E^c(z)\oplus \RR X(z))$ for
$z\in\Lambda$ with the exception of the flow direction
(along which the Lyapunov exponent is zero).

In this way we show that the attractor for the flow of
$X$ is chaotic, in the sense that it admits a physical
probability measure with $k+1$ positive Lyapunov exponents.

\emph{This concludes the proof of Theorem~\ref{mthm:pri}}.

\subsection{Unfolding $X$.}
\label{sec:unfold-x}

The $E^{ss}$ direction on $\Sigma_2$ can be made uniformly
contracting with arbitrarily strong contraction rate (see
Section~\ref{sec:unpert-singul-flow} and
Remark~\ref{rmk:domination}). Moreover, the $x_2$
direction is dominated by all the other
directions under the flow $X^t$. Thus $E^{ss}$ 
is a stable direction for the flow over $\Lambda$ which is
dominated by any complementary direction.

Hence this uniformly contracting foliation $\cF^{ss}_X$ admits
a continuation $\cF^{ss}_Y$ for all flows $Y^t$ where $Y$ is
close to $X$ in $\fX^2(M)$, that is, in the $C^2$
topology. Since it is a one-dimensional foliation whose
contraction rate can be made arbitrarily small, the
holonomies along its leaves are of class $C^{\gamma}$ for
some $\gamma>1$, see \cite[Theorem 6.2]{HP70}.

This ensures that the return map $R_Y$ to the cross-section
$I_\epsilon\times I_\epsilon\times\{1\}\times B^{k+2}$,
where $I_\epsilon:=(-1-\epsilon,1+\epsilon)$, admits a
one-dimensional invariant foliation such that the quotient
map $g_Y$ of $R_Y$ on the leaves of this foliation is a
$(k+3)$-dimensional $C^{\gamma}$ map, for \emph{any vector
  field $Y$ close to $X$ in the $C^2$ topology}. In
addition, the leaves of $\cF^{ss}_Y$ are $C^1$-close to
the leaves of the original $\cF^{ss}_X$ foliation.

Considering the set $I_\epsilon\times\{(0,1)\}\times
B^{k+2}$ diffeomorphic to $I_\epsilon\times B^{k+2}$, we see
that $\cF^{ss}_X$ is transverse to this set and thus the
continuation remains transverse. Hence we can see the map
$g_Y$ as a map between subsets of $I_\epsilon\times
B^{k+2}$. Let $\pi^{ss}_Y:\Sigma_2\to I_\epsilon\times
B^{k+2}$ be the projection along the leaves of $\cF^{ss}_Y$
in what follows and let $\ell:=\Sigma_2\cap
W^s_{loc}(s_0(Y))$ be the connected component of the local
stable manifold of the continuation of $s_0$ for $Y$ on the
cross-section $\Sigma_2$.

We have scarce information about the dynamics of this map:
it has a sink $p(Y)$ (the continuation of the sink of $R$)
and is $C^\gamma$ close to the quotient of the map $R$ over
the foliation $\cF^{ss}_X$. Thus, for points outside a
neighborhood of $\pi_Y^{ss}(\ell)$ and away from the stable
manifold of $p(Y)$, we should have ``hyperbolicity'' for
$g_Y$ due to proximity of $R_Y$ to $R$, that is, there is a
pair of complementary directions on the tangent space such
that one is expanded and the other contracted by the
derivative of $g_Y$. The interplay of this hyperbolic-like
behavior with the behavior near $\ell$ is unknown to us.

\begin{conjecture}\label{conj:physicalViana}
  Similarly to the one-dimensional setting, \emph{$g_Y$
    admits a physical hyperbolic measure $\mu_Y$ for all
    those vector fields $Y$ which are $C^2$ close to $X$ and
    the stable manifold of the sink does not contain the
    critical region.} Moreover the basin of $\mu_Y$ should
  be the complement of the stable manifold of the sink.
\end{conjecture}

However, we can be more specific along certain submanifolds
of the space of vector fields, as follows.

\subsubsection{Keeping the domination on $\Sigma_2$ under
  perturbation}
\label{sec:keeping-dominat}

The argument presented in the
Section~\ref{sec:dominat-splitt}, proving smoothness of a
$3$-dimensional stable foliation $\cF$ after quotienting by
$\pi_1$, strongly depends on the fact that $\pi_1$
identifies every point whose first coordinate is $1$, which
is represented by an infinite contraction there. We just
have to consider \eqref{eq:crucial}, which holds because the
point $0$, corresponding to the intersection of $\Sigma_2$
with the local stable manifold of $s_0$, is sent by $g_0$ to
$1$ on each side, that is
$g_0(0^\pm)=\lim_{t\to0^\pm}g_0(t)=1$.

In order to keep the strong domination for a perturbation
$Y$ of $X$ in the $C^2$ topology, \emph{we restrict the
  perturbation in such a way that the corresponding points
  $P_2(Y)$ and $P_3(Y)$ are in the same stable leaf
  $\xi\in\cF_Y$}. This is well defined according to
Remark~\ref{rmk:1dfoliation}, see Figure~\ref{fig:startflow}
for the positions of $P_2$ and $P_3$.  Here $P_2(Y)$ and
$P_3(Y)$ are the points of first intersection of each branch
$W^u(s_0)\setminus\{s_0\}$ of the one-dimensional unstable
manifold of the equilibrium $s_0$. We note that we can write
each branch as an orbit of the flow of $Y$ and so the notion
of first intersection with $\Sigma_2$ is well defined. This
restriction on the vector field corresponds to restricting
to a $(k+4)$-codimension submanifold $\cP$ of the space of
vector fields $\fX^2(M)$.

In this way, on the one hand, the same arguments of
Section~\ref{sec:dominat-splitt} can be carried through and
the strong domination persists for vector fields
$Y\in\cP$. On the other hand, this implies that there exists
a $R_Y$-invariant $3$-dimensional contracting $C^\gamma$
foliation $\cF_Y$ of $\Sigma_2$, for some $\gamma>1$, with
$C^1$ leaves, for all vector fields $Y$ close enough to $X$
within $\cP$.

We can then quotient $R_Y$ over the leaves of $\cF_Y$ to
obtain a $(k+1)$-dimensional map $g_Y$.  We note that,
defining the cylinder
$\cC:=I_\epsilon\times\{0\}\times\{1\}\times
e(\TT^k\times0)$ (diffeomorphic to $I_\epsilon\times\TT^k$)
inside $I_\epsilon\times I\times\{1\}\times B^{k+2}$ (recall
that $I_\epsilon=(-1-\epsilon,1+\epsilon)$) we have that the
initial foliation $\cF$ is everywhere transverse to $\cC$.
Therefore, since $\cC$ is a proper submanifold, the
continuation $\cF_Y$ is still transverse to $\cC$ for
$Y\in\cN\cap\cP$. Hence we can define a corresponding
quotient map $g_Y:I_\epsilon\times\TT^k\circlearrowleft$
which will \emph{not be}, in general, either a direct or
skew-product along the $I_\epsilon$ and $\TT^k$ directions.

We observe that $g_Y$ is close to $f_0\times f_1$ on
$I\times\TT^k$, recall Remark~\ref{rmk:cylinder} during the
construction of the original flow. Hence for pieces of
orbits which remain away from a neighborhood of the critical
set and away from the basin of the sink, we have uniform
expansion in all directions (akin to condition C on the
statement of Theorem~\ref{thm:Misiurew}). 

% \begin{remark}\label{rmk:C1+quotient}
%   We stress that since we are in the setting of strong
%   domination, the foliation $\cF$ has $C^\gamma$ holonomies
%   for some $\gamma>1$. Thus the quotient $g_Y$ is a $C^{\gamma}$
%   map.
% \end{remark}

\subsubsection{Keeping the saddle-connection}
\label{sec:nearby-vector-fields}

In addition to keeping the foliation $\cF$, we may impose
the restriction already mentioned in
Remark~\ref{rmk:keepsaddleconnection} to \emph{keep also the
  connection between $s_0$ and $s_1$}: the component of the
unstable manifold of $s_0$ through $\Sigma^+$ (i.e. the
orbit of $P_2(Y)$) is contained in the stable manifold of
$s_1$, a $(k+1)$-codimension condition on the family of all
vector fields $Y$ $C^2$ close to $X$. Let $\cN$ be the
submanifold of such vector fields in a neighborhood of $X$.

Therefore we can ensure that there exists a stable foliation
$\cF_Y$ nearby $\cF_X$ for every vector field
$Y\in\cN\cap\cP$, invariant under the corresponding return
map $R_Y$. 
We can again quotient $R_Y$ over the leaves of
$\cF_Y$ to obtain a $(k+1)$-dimensional map $g_Y$. We
observe that $\cN\cap\cP$ will have codimension $2k+5$ since
the conditions defining $\cN$ and $\cP$ are independent.

% We can argue as in Rovella~\cite{Ro93} further restricting
% the vector field $Y$ to a codimension $2$ submanifold
% \begin{align*}
%   \cN_0:=\left\{ Y\in\cN : 1\pm \text{  is preperiodic
%       repelling under the quotient map  } g_Y
%     \right\}
% \end{align*}
%%%%%%%%%%%%%%%%%%%%%%%%%%%%%%%%%%%%%%%%%%%%%%%%%%%%%%%

\subsubsection{Keeping the one-dimensional quotient map}
\label{sec:symmetr-unfold-x}

We can also perturb the vector field $X$ within the manifold
$\cN\cap\cP$ keeping the saddle-connection in such a way
that we obtain a one-dimensional $C^{1+}$ quotient map. In
this setting we can apply Benedicks-Carleson exclusion of
parameters techniques along these families of flows, exactly
in the same way Rovella proved his main theorem in
\cite{Ro93}. \emph{So we have an analogous result to
  Rovella's if we perturb the flow keeping the symmetry
  which allows us to project to a one-dimensional map, that
  is, if $g_Y$ is a skew-product over $I_\epsilon$.}

\begin{figure}[htpb]
  \centering
  \psfrag{0}{$0$}\psfrag{p}{$p_1$}\psfrag{r}{$p_0=-1$}
\psfrag{s1}{$s_1$}\psfrag{P0}{$\pi_2$}
\psfrag{N}{$N$}\psfrag{n}{$+1$}\psfrag{s}{$-1$}
\psfrag{T}{$\TT_1^k$}\psfrag{sr}{$\hat s$}
  \includegraphics[scale=0.5]{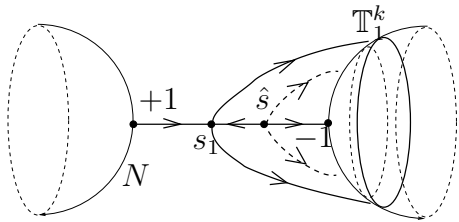}
  \caption{\label{fig:symperturb}A symmetrical unfolding. We
    can still project the dynamics into a one-dimensional
    map.}
\end{figure}

Since we have a well defined strong
stable foliation for the Poincar\'e map, we can quotient
out along such stable foliation obtaining a map in the
torusphere. By {\em keeping the symmetry} we mean that we
unfold our Misiurewicz type flow preserving the invariance
of each parallel torus in the torusphere.  We can consider
families which unfold the criticality introduced by the
singularity $s_1$ (as in the right hand side of
Figure~\ref{fig:contractingLorenz}) and/or unfold the
intersection of the unstable manifold of $s_1$ with
$\TT_1^k$, see Figure~\ref{fig:symperturb}.  Once more, this
permits us to reduce the study of the attractor to a
one-dimensional problem.

In the very small manifold within $\cP\cap\cN$ where the
one-dimensional quotient map is kept, we can argue just like
Rovella in \cite{Ro93} obtaining the same results.

%%%%%%%%%%%%%%%%%%%%%%%%%%%%%%%%%%%%%%%%%%%%%%%%%%%%

\subsubsection{Loosing the one-dimensional quotient map}
\label{sec:non-symmetr-unfold}

\begin{figure}[htpb]
  \centering
  \psfrag{0}{$0$}\psfrag{p}{$p_1$}\psfrag{r}{$p_0=-1$}
\psfrag{s1}{$s_1$}\psfrag{P0}{$\pi_2$}
\psfrag{N}{$N$}\psfrag{n}{$+1$}\psfrag{s}{$-1$}
\psfrag{T}{$\TT_1^k$}\psfrag{sr}{$\hat s$}
  \includegraphics[scale=0.5]{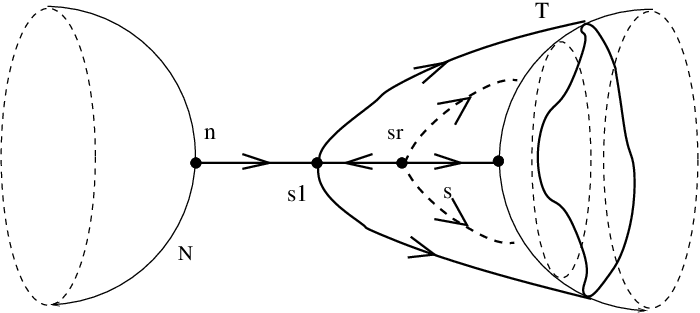}
  \hspace{1cm}
  \includegraphics[scale=0.5]{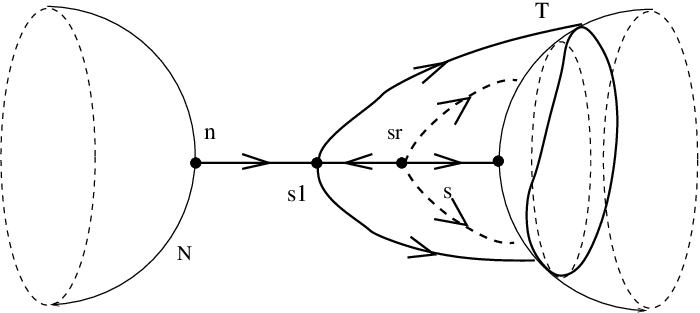}
  \caption{\label{fig:skewperturb} On the left: a
    non-symmetrical perturbation. We cannot reduce the
    analysis to a one-dimensional map. On the right: An open
    region nearby the singularity is sent into the basin of
    a sink. }
\end{figure}

We note also that, even if we keep the saddle-connection, it
is very easy to \emph{perturb this flow} to another
arbitrarily close flow such that \emph{the quotient to a
  one-dimensional map is not defined}, as the left hand side
of Figure~\ref{fig:skewperturb} suggests.

\subsubsection{Breaking the saddle-connection}
\label{sec:breaking-saddle-conn}

Another way to break the one-dimensional quotient is to
break the saddle-connection. Then the orbit of $P_2(Y)$ is
no longer contained in the stable manifold of $s_1$. So the
positive orbit of $P_2(Y)$ under $Y^t$ follows the unstable
manifold of $s_1$ and crosses $\Sigma_2$ at some point, but
points in the image $T^+_Y(\Sigma^+)$ near $P_2(Y)$ may no
longer be mapped preserving the stable foliation of the
solenoids (see Figure~\ref{fig:startflow}).
In this case we still have a quotient map
$g_Y:I_\epsilon\times\TT^k\circlearrowleft$ but it is no longer a
skew-product over $I_\epsilon$.

\subsubsection{Returning away from the sink}
\label{sec:return-away-from}

We conjecture that, in the situation depicted in the left
side of Figure~\ref{fig:skewperturb}, that is, we have a
non-skew-product map $g_Y$ but the image of $P_2(Y)$ under
$R_Y$ does not intersect the basin of the sink $P_0$,
\emph{either with or without a saddle-connection} between
$s_0$ and $s_1$, then it is \emph{always true that there
  exists an expanding (all Lyapunov exponents are positive)
  physical measure with full basin outside the basin of the
  sink}.

The motivation behind this conjecture is that $g_Y$ is close
to $g_0$ and away from the singular set and away from the
basin of the sink we have uniform expansion, since this
behavior was present in the original map $g_0$ and is
persistent. Therefore, we have the interplay between
expansion and a critical region which will be approximately
a circle on the cylinder $I_\epsilon\times \TT^k$. This
setting is similar to the one introduced by Viana in
\cite[Theorem A]{Vi97} and carefully studied in
e.g. \cite{alves-viana2002,AA03,alves-luzzatto-pinheiro2005}.

\subsubsection{Returning inside the basin of the sink}
\label{sec:return-inside-basin}

On the one hand, the perturbation can send points nearby
$P_2$ into the basin of the sink $P_0$, \emph{either with or
  without a saddle-connection} between $s_0$ and $s_1$, as
shown in the right hand side of
Figure~\ref{fig:skewperturb}. We prove, in
Section~\ref{sec:full-basin-attract}, that this implies that
the \emph{basin of the sink grows to fill the whole
  manifold} Lebesgue modulo zero, for such vector fields in
the submanifold $\cP\cap\cN$.

% \begin{figure}[htpb]
%   \centering
%   \psfrag{0}{$0$}\psfrag{p}{$p_1$}\psfrag{r}{$p_0=-1$}
% \psfrag{s1}{$s_1$}\psfrag{P0}{$\pi_2$}
% \psfrag{N}{$N$}\psfrag{n}{$+1$}\psfrag{s}{$-1$}
% \psfrag{T}{$\TT_1^k$}\psfrag{sr}{$\hat s$}
%   \includegraphics[scale=0.5]{tubularskewperturb3.eps}
%   \caption{\label{fig:sinkperturb} An open region nearby the
%     singularity is sent into the basin of a sink.}
% \end{figure}

\emph{This completes the proof of
  Theorem~\ref{mthm:codimension}} except for the last item
proved in Section~\ref{sec:full-basin-attract}.

%%%%%%%%%%%%%%%%%%%%%%%%%%%%%%%%%%%%%%%%%%%%%%%%%%%%%%%%

\section{Higher dimensional Benedicks-Carleson conditions}
\label{sec:higher-dimens-misiur}

To present the statement and proof of the existence of an
absolutely continuous invariant probability measure on the
setting of higher dimensional maps, we need to recall some
notions from non-uniformly expanding dynamics.

\subsection{Non-flat critical or singular sets}
\label{sec:non-flat-critic}

Let $f$ be a $C^{1+}$ local diffeomorphism outside a compact
proper submanifold $\cS$ of $M$ with positive codimension.
 The set $\cS$ may be taken as some set of critical points
 of $f$ or a set where $f$ fails to be differentiable. The
 submanifold $\cS$ has an at most countable collection of
 connected components $\{\cS_i\}_{i\in\NN}$, which may have
 different codimensions. This is enough to ensure that the
 volume or Lebesgue measure of $\cS$ is zero and, in
 particular, that $f$ is a \emph{regular map}, that is, if
 $Z\subset M$ has zero volume, then $f^{-1}(Z)$ has zero
 volume also.

 \emph{In what follows we assume that the number of
   connected components of $\cS$ is finite.} It should be
 possible to drop this condition if we impose some global
 restrictions on the behavior of the map $f$, see
 \cite{PRV98} for examples with one-dimensional ambient
 manifold $M$ of non-uniformly expanding maps with
 infinitely many critical points. We do not pursue this
 issue in this work.

 We say that \( \cS \subset M\) is a {\em non-flat critical
   or singular set} for $f$ if the following conditions hold
 for each connected component $\cS_i$. The first one
 essentially says that $f$ {\em behaves like a power of the
   distance} to \( \cS_i \): there are constants $B>1$ and
 real numbers $\alpha_i>\beta_i$ such that
 $\alpha_i-\beta_i<1, 1+\beta_i>0$ and \emph{on a
   neighborhood $U_i$ of $\cS_i$} (where
 $U_i\cap\cS_j=\emptyset$ if $j\neq i$ and we also take
 $U_i\subset B(\cS_i,1/2)$) for every $x\in U_i$
\begin{enumerate}
 \item[(S1)]
   $\displaystyle{\frac{1}{B}\dist(x,\cS)^{\alpha_i}\leq
     \frac {\|Df(x)v\|}{\|v\|}\leq B\dist(x,\cS)^{\beta_i}}$
   for all $v\in T_x M$.
\end{enumerate}
Moreover, we assume that the functions \(  \log|\det Df(x)| \) and
\( \log \|Df(x)^{-1}\| \) are \emph{locally Lipschitz} at points
\( x\in U_i \) with Lipschitz constant
depending on $\dist (x, \cS)$: for every $x,y\in
U_i$ with $\dist(x,y)<\dist(x,\cS_i)/2$ we have
\begin{enumerate}
\item[(S2)] $\displaystyle{\left|\log\|Df(x)^{-1}\|-
\log\|Df(y)^{-1}\|\:\right|\leq
\frac{B}{\dist(x,\cS_i)^{\alpha_i}}\dist(x,y)}$;
\item[(S3)] $\displaystyle{\left|\log|\det Df(x)|- \log|\det
      Df(y)|\:\right|\leq
    \frac{B}{\dist(x,\cS_i)^{\alpha_i}}\dist(x,y)}$.
 \end{enumerate}
 The assumption that the number of connected components is
 finite implies that there exists $\beta>0$ such that
 $\max_i\{|\alpha_i|,|\beta_i|\}\le\beta$ and we also assume
 that for all $x\in M\setminus\cS$
 \begin{itemize}
 \item[(S4)] $\displaystyle{\frac{1}{B}\dist(x,\cS)^{\beta}\leq
     \frac {\|Df(x)v\|}{\|v\|}\leq B\dist(x,\cS)^{-\beta}}$
   for all $v\in T_x M$.
 \end{itemize}
 The case where $\cS$ is equal to the empty set may also be
 considered. The assumption $1+\beta_i>0$ prevents that the
 image of arbitrary small neighborhoods around the singular
 set accumulates every point of $M$ (consider e.g. the Gauss
 map $[0,1]\ni x\mapsto x^{-1}\bmod1$) which would prevent a
 meaningful definition of ``singular value'', see the
 statement of Theorem~\ref{thm:Misiurew} in what follows.

\subsection{Hyperbolic times}
\label{sec:hyperb-times}

We write $x_i=f^i(x)$ for $x\in M$, $i\ge0$, and also
$S_n^f\varphi(x)=S_n\varphi(x)=\sum_{i=0}^{n-1}
\varphi(x_i), n\ge1$ for the ergodic sums of a function
$\varphi:M\to\RR$ with respect to the action of $f$, in what
follows.  For the next definition it will be useful to
introduce \( \dist_{\delta}(x,\cS) \), the \( \delta
\)-\emph{truncated} distance from \( x \) to~$\cS$, defined
as \( \dist_{\delta}(x,\cS) = \dist(x,\cS) \) if \(
\dist(x,\cS) \leq \delta\), and \( \dist_{\delta}(x,\cS) =1
\) otherwise.  From now on we write
$\psi(x)=\log\|Df(x)^{-1}\|$ and $\fD_r(x)=-\log
\dist_r(x,\cS)$ for $x\in M\setminus\cS$ and $r>0$.

Let $\beta>0$ be given by the non-flat conditions on $\cS$,
and fix $b>0$ such that $b < \min\{1/2,1/(4\beta)\}$.  Given
$c>0$ and $\delta>0$, we say that $h$ is a {\em
  $(e^c,\delta)$-hyperbolic time}
% \footnote{In the case
%     $\cs=\emptyset$ the definition of
%     $(\sigma,\delta)$-hyperbolic time reduces to the first
%     condition in \eqref{d.ht}, and we simply call it
%     $\sigma$-hyperbolic time. }}
for a point $x\in M$ if, for all $0\le k \le h$,
 \begin{equation}\label{d.ht}
%\prod_{j=h-k}^{h-1}\|Df(x_j))^{-1}\| \le \sigma^k
S_k\psi(x_{n-k})\le -ck
\quad\text{and}\quad
% \prod_{j=n-k}^{n}
S_k\fD_\delta(x_{n-k})\le b c k
%\dist_\delta(x_{h-k}, \cS)\ge \sigma^{b k}.
 \end{equation}
 We convention that an empty sum evaluates to $0$ so
 that the above inequalities make sense for all indexes in
 the given range.

 We say that the {\em frequency of
   $(e^c,\delta)$-hyperbolic times} for $x\in M$ is
 greater than $\theta>0$ if, for infinitely many times $n$,
 there are $h_1<h_2\dots <h_\ell\le n$ which are
 $(e^c,\delta)$--hyperbolic times for $x$ and $\ell
 \ge\theta n$.

 The following statement summarizes the main properties of
 hyperbolic times. For a proof the reader can consult
 \cite[Lemma 5.2]{ABV00} and \cite[Corollary 5.3]{ABV00}.

 \begin{proposition}\label{p.contr}
   There are $0<\delta_1<\delta/4$ and $C_1>0$ (depending
   only on $\delta$ and $\sigma$) such that if $h$ is a
   $(e^c,\delta)$-hyperbolic time for $x$, then there is
   a \emph{hyperbolic neighborhood} $V_x$ of $x=x_0$ in $M$
   for which
   \begin{enumerate}
   \item $f^h$ maps $V_x$ diffeomorphically onto the ball of
     radius $\delta_1$ around $x_n$;
   \item for $1\le k \le h$ and $y, z\in V_x$, $
     \dist(y_{n-k},z_{n-k}) \le
     C_1 e^{ck/2}\dist(y_h,z_h)$;
   \item for all $y\in V_x$, $S_h\fD_\delta(y)\le
     S_h\fD_\delta(x) + o(\delta)$ where
     $o(\delta)/\delta\to0$ when $\delta\to0$;
   \item $f^h\vert V_x$ has distortion bounded by $C_1$: if
     $y, z\in V_x$, then
     $ \frac{|\det Df^h (y)|}{|\det Df^h (z)|}\le C_1. $
\end{enumerate}
\end{proposition}

\begin{remark}
  \label{rmk:preball-away}
  The image of $V_x$ by $f^h$ is away from
  $B(\cS,3\delta/4)$: for $y\in V_x$
  \begin{align*}
    \dist(y_h,\cS)
    \ge
    \dist(x_h,\cS) - \dist(y_h,x_h)
    \ge \delta-\delta_1>\frac34 \delta
  \end{align*}
  since $\dist_\delta(x_h,\cS)\ge1$.
\end{remark}

Item (3) above, which is not found in \cite{ABV00}, is an
easy consequence of item (2): for every $1\le k\le h$ and
$y\in V_x$
\begin{align*}
  -\log\frac{\dist_\delta(y_{h-k},\cS)}
  {\dist_\delta(x_{h-k},\cS)} &\le
  -\log\frac{\dist_\delta(x_{h-k},\cS)-
    \dist(y_{h-k},x_{h-k})}
  {\dist_\delta(x_{h-k},\cS)}
  \\
  &= -\log\left( 1-\frac{\dist(y_{h-k},x_{h-k})}
    {\dist_\delta(x_{h-k},\cS)} \right)
  \\
  &\le -\log\left(1-\frac{\delta_1\sigma^{k/2}}
    {e^{bck}} \right) \le
  -\log\left(1-\frac34\delta e^{(1/2-b)ck}\right)
  \\
  &\le \frac{(3/4)\delta}{1-(3/4)\delta} \cdot \frac34\delta
  e^{(1/2-b)ck}.
  \end{align*}
  So each point $y$ in a hyperbolic neighborhood as above
  satisfies
  \begin{align*}
    S_h\fD_\delta(y)
    &\le
     S_h\fD_\delta(x)
    +
    \frac{(3/4)^2\delta^2}{1-(3/4)\delta}
    \sum_{j=0}^h e^{(1/2-b)cj}
    \le
    S_h\fD_\delta(x)
    +
    \frac{(3/4)^2\delta^2}{1-(3/4)\delta}\cdot
    \frac1{1-e^{bc}},
  \end{align*}
  and the last term is $o(\delta)$.

In what follows we say that an open set $V$ is a
\emph{hyperbolic neighborhood} of any one of its points
$x\in V$, \emph{with $(e^{c/2},\delta)$-distortion time
  $h$} (or, for short, a
$(e^{c/2},\delta)$-\emph{hyperbolic neighborhood}), if
the properties stated in items (1) through (4) of
Proposition~\ref{p.contr} are true for $x$ and $k=h$.

We also say that a given point $x$ has \emph{positive
  frequency of hyperbolic neighborhoods} bounded by
$\theta>0$ if there exist $c,\delta>0$ and
neighborhoods $V_{h_k}$ of $x$ with
$(e^{c},\delta)$-distortion time $h_k$, for all $k\ge1$,
such that for every big enough $k$ we have $k\ge\theta h_k$.

\subsection{Existence of many hyperbolic neighborhoods
  versus absolutely continuous invariant probability
  measures}
\label{sec:existence-hyperb-tim}

Hyperbolic times appear naturally when $f$ is assumed to be
{\em non-uniformly expanding} in some set $H\subset M$:
there is some $c>0$ such that for every $x\in H$ one has
 \begin{equation} \label{liminf1}
\liminf_{n\to +\infty}\frac{1}{n}
S_n\psi(x)<-c,
 \end{equation}
 and points in $H$ satisfy some {\em slow recurrence to the
   critical or singular set}: given any $\epsilon>0$ there
 exists $\delta>0$ such that for every $x\in H$
 \begin{equation}\label{eq:slow-recurrence}
   \limsup_{n\to+\infty} \frac1{n} S_n\fD_{\delta}(x)
   \le\epsilon.
\end{equation}

The next result has been proved in \cite[Theorem~C \& Lemma
5.4]{ABV00}. It provides sufficient conditions for the
existence of many hyperbolic times along the orbit of points
satisfying the non-uniformly expanding and slow recurrence
conditions.

\begin{theorem}\label{thm:abv0}
  Let $f: M\to M$ be a $C^{1+}$ local diffeomorphism outside
  a non-flat critical or singular set $\cS\subset M$. If
  there is some set $H\subset M\setminus\cS$ such
  that~\eqref{liminf1} and~\eqref{eq:slow-recurrence} hold
  for all $x\in H$, then for any given $0<\xi<1$ and
  $0<\zeta<bc$ there exist $\delta>0$ and $\theta>0$ such
  that the frequency of $(e^{-c\xi},\delta)$-hyperbolic
  times for each point $x\in H$ is bigger than
  $\theta$. Moreover for such hyperbolic times $h$ we have
  $S_h\fD_\delta(x) \le \zeta h.$
\end{theorem}

Together with Proposition~\ref{p.contr} the results from
Theorem~\ref{thm:abv0} ensure the existence of positive
frequency of hyperbolic neighborhoods around Lebesgue almost
every point.

This will imply the existence of absolutely continuous
invariant probability measures for the map $f$ through the
following result from~\cite{Pinheiro05}.

\begin{theorem}
  \label{thm:abv1}
  Let $f: M\to M$ be a $C^{1+}$ local diffeomorphism outside
  a non-degenerate exceptional set $\cS\subset M$. If there
  are $c,\delta>0$
% , and $H\subset M\setminus\cS$
%   with $m(H)>0$
  such that the frequency of $(e^c,\delta)$-hyperbolic
  neighborhoods is bigger than $\theta>0$ for Lebesgue
  almost every $x\in M$, then $f$ has some absolutely
  continuous invariant probability measure.
\end{theorem}

Since we are using a modified definition of hyperbolic time,
we present a proof of Theorems~\ref{thm:abv0}
and~\ref{thm:abv1} in Section~\ref{sec:hyperb-neighb-constr}
for completeness.

\subsection{Existence of absolutely continuous probability
  measures}
\label{sec:existence-absolut-co}

The following theorem provides higher dimensional existence
result for physical measures which applies to our
setting.

\begin{theorem}
  \label{thm:Misiurew}
  Let $f:M\setminus\cS \to M$ be a $C^{1+}$ local
  diffeomorphism, where $\cS$ is a compact sub-manifold of
  $M$ which is a non-flat critical/singular set for $f$.  We
  define
  $$f(\cS):=\cap_{n\ge1}\overline{f\big(B_{1/n}(\cS)\big)},$$
  the set of all accumulation points of sequences $f(x_n)$
  for $x_n$ converging to $\cS$ as $n\to+\infty$, and assume
  that $f$ also satisfies:
  \begin{itemize}
  \item[A:] \emph{$f$ is non-uniformly expanding along the
      orbits of critical values}: there exist $c_0>0$ and
    $N\ge1$ such that for all $x\in f(\cS)$ and $n\ge N$ we
    have $S_n\psi(x) \le -c_0 n$;
  \item[B:] \emph{the critical set has slow recurrence to
      itself}: given $\epsilon>0$ we can find $\delta>0$
    such that for all $x\in f(\cS)$ there exists $N=N(x)$
    satisfying $S_n\fD_\delta(x)\le\epsilon n$ for
        every $n\ge N$;
  \item[C:] \emph{$f$ is uniformly expanding away from the
      critical/singular set}: for every neighborhood $U$ of
    $\cS$ there exist $c=c(U)>0$ and $K=K(U)>0$
    such that for any $x\in M$ and $n\ge1$ satisfying
      $x=x_0, x_1,\dots, x_{n-1}\in M\setminus U$,
      then $S_n \psi(x) \le K- c n$.
    \item[D:] \emph{$f$ does not contract too much when
        returning near the critical/singular set}: that is,
      there exists $\kappa>0$ and a neighborhood $\hat U$ of
      $\cS$ such that for every open neighborhood $U\subset
      \hat U$ of $\cS$ and for $x=x_0$ satisfying $x_1,
      \dots, x_{n-1}\in M\setminus U$ and either
      $x_0\in U$ or $x_n \in  U$, then
      $S_n\psi(x_0)\le\kappa$.
  \end{itemize}
  Then $f$ has an absolutely continuous invariant
  probability measure $\nu$ such that $\fD_d$ is
  $\nu$-integrable for some (and thus all) $d>0$.
\end{theorem}

We observe that condition B above ensures, in particular,
that $f(\cS)\cap\cS=\emptyset$, for otherwise
$\fD_\delta(x)$ is not defined for $x\in\cS\cap
f(\cS)$. Moreover, condition C above is just a convenient
translation to this higher dimensional setting of the
conclusion of the one-dimensional theorem of Ma\~n\'e
\cite{Man85}, ensuring uniform expansion away from the
critical set and basins of periodic attractors. It can be
read alternatively as: given $\delta>0$ there are
$C,\lambda>0$ such that if $x_i\in M\setminus B(\cS,\delta)$
for $i=0,\dots,n-1$, then $\|Df^n(x)^{-1}\|\le C e^{-\lambda
  n}$.  In addition, condition D is a translation to our
setting of a similar property that holds for unidimensional
multimodal ``Misiurewicz maps'', that is, for maps whose
critical orbits are non-recurrent, ensuring a minimal lower
bound for the derivative of the map along orbits which
return near $\cS$.

Now we show that under the conditions in the statement of
Theorem~\ref{thm:Misiurew}, we can find a full measure
subset of points of $M$ having positive density of
hyperbolic times.

\begin{theorem}
  \label{thm:dense-hyp-times}
  Let $f:M\setminus\cS \to M$ be a $C^{1+}$ local
  diffeomorphism away from a non-flat critical/singular set
  $\cS$, satisfying all conditions in the statement of
  Theorem~\ref{thm:Misiurew}. Then for every small enough
  $0<\xi<1$ there exists $\delta=\delta(\xi)>0$ and
  $\theta=\theta(\xi,\delta)>0$ such that Lebesgue almost
  every $x\in M$ admits positive frequency bounded by
  $\theta$ of $(e^{-\xi c_0},\delta)$-hyperbolic
  neighborhoods.
\end{theorem}

From this result we deduce Theorem~\ref{thm:Misiurew}
appling Theorem~\ref{thm:abv1}. So all we need to do is
prove Theorem~\ref{thm:dense-hyp-times}. For the
integrability of $\fD$ see
Remark~\ref{rmk:log-dist-integrable} in what follows.

\subsubsection{Existence of hyperbolic neighborhoods}
\label{sec:existence-hyperb-nei}

Fix $\xi_0,\epsilon,\tilde\delta>0$ and small enough so that
condition B is satisfied in what follows and $\xi_0c_0<b$.
Let $\zeta>0$ be small enough in order that
  \begin{align}\label{eq:zeta}
    \frac{\xi_0 c_0}{1+\alpha_i} + 2\zeta
    < \frac{\xi_0 c_0}{1+\zeta} \quad\text{for all   } i.
  \end{align}
  Depending on $f$ and $\xi_0$, we can choose the pair
  $(\epsilon,\tilde\delta)$ so that, from conditions A and B
  above together with Theorem~\ref{thm:abv0}, every point
  $z\in f(\cS)$ has infinitely many $(e^{-2\xi_0
    c_0},\tilde\delta)$-hyperbolic times $h_1<h_2<h_3<\dots$
  satisfying
  \begin{align}\label{eq:Daverage}
    S_{h_i}\fD_{\tilde\delta}(z)\le\zeta h_i
    \quad\text{for}\quad i\ge1.
  \end{align}
  Then there are corresponding hyperbolic neighborhoods
  $V_i$ of $z$ satisfying the conclusions of
  Proposition~\ref{p.contr} for each hyperbolic time $h_i$
  of $z$, where $\tilde\delta>0$ does not depend on
  $z\in\cS$, that is $f^{h_i}\mid V_i: V_i\to
  B(z_{h_i},\delta_1)$ is a diffeomorphism with a ball of
  radius $\delta_1 \in (0,\tilde\delta/4)$ whose inverse is
  a contraction with rate bounded by $e^{-\xi_0 c_0 h_i}$.

  We will consider, instead of $V_i$, the subset $B_i\subset
  V_i$ given by
  \begin{align}
    \label{eq:halfball}
    B_i = \big(f^{n_i}\mid
    V_i\big)^{-1}(B(z_n,\delta_2))
  \end{align}
  where we set $2\delta_2=\delta_1$.  Since $f(\cS)$ is
  compact, we can cover this set by finitely many hyperbolic
  neighborhoods of the type $B_i$.

  Now we fix a connected component $\cS_i$ of $\cS$.

  \subsubsection{Hyperbolic neighborhoods near the
    critical/singular set}
\label{sec:hyperb-neighb-near}

Let $T_i$ be the smallest distortion (or hyperbolic) time
associated to the balls covering $f(\cS_i)$.  We remark that
$T_i$ can be taken arbitrarily big, independently of
$\xi_0,\zeta,\tilde\delta$, because every $z\in f(\cS_i)$
has positive frequency of hyperbolic neighborhoods and,
consequently, the open neighborhoods in the above covering
can be made arbitrarily small.

We observe that, by definition of hyperbolic neighborhoods,
if $\dist(z,f(\cS_i))<e^{-b j \xi_0 c_0}$ for some $j\ge1$
and there exists $x\in f(\cS_i)$ such that $z\in V_i(x)$,
then the corresponding distortion time $h_i$ satisfies
$h_i\ge j$.

We note also that, by the non-flat condition (S1) on $f$
near $\cS_i$, a $\rho$-neighborhood of $\cS$ is sent into a
$\rho^{1+\beta_i}$-neighborhood of $f(\cS_i)$, for each
$\rho>0$. Indeed, since $\cS_i$ is assumed to be a
submanifold, we can find for $x$, on a tubular neighborhood
of $\cS_i$ with radius $\rho$, a curve $\gamma:[0,\rho]\to
M$ from $\gamma(0)\in\cS_i$ to $\gamma(1)=x$ such that
$\dist(\gamma(t),\cS_i)=t$ and $\|\dot\gamma(t)\|=1$ for
$t\in[0,\rho]$. Hence
  \begin{align}
    \dist(f(x),f(\cS_i))
    &\le\nonumber
    \int_0^\rho
    \|Df(\gamma(t))\dot\gamma(t)\|\,dt
    \le
    \int_0^\rho B \dist(\gamma(t),\cS)^{\beta_i}\,dt
    \\
    &=
    \frac{B\rho^{1+\beta_i}}{1+\beta_i}
    =\frac{B}{1+\beta_i}\dist(x,\cS)^{1+\beta_i}.
    \label{eq:distS1}
  \end{align}
  Therefore we can find $C_2=C_2(\delta_2)>0$ such that
  $C_2^{1+\beta_i} B/(1+\beta_i)=\delta_2$ and if, for some
  $j\ge0$
  \begin{align}\label{eq:dist-time}
    \dist(x,\cS_i)\le C_2 e^{-\xi_0 c_0
      (T_i+j)/(1+\alpha_i)}:=d_j,
  \end{align}
  then the smallest possible distortion time $h$ of $f(x)$
  is at least $\frac{1+\beta_i}{1+\alpha_i}(T_i+j)$, since
  \begin{align*}
    \delta_2 e^{-\xi_0 c_0 h}
    \le
    \dist(f(x),f(\cS_i))
    \le
    \delta_2 e^{- \xi_0 c_0
      \frac{1+\beta_i}{1+\alpha_i}(T_i+j)}
    \le
    \delta_2 e^{-b \xi_0 c_0
      \frac{1+\beta_i}{1+\alpha_i}(T_i+j)}.
  \end{align*}
  Note that by the previous observations and
  Remark~\ref{rmk:preball-away} we have
  \begin{align*}
    f(B(\cS_i,\tilde\delta)\subset
    B\left(f(\cS_i),\frac{B}{1+\beta_i}
      \tilde\delta^{1+\beta_i}
    \right)
  \end{align*}
  and we can assume without loss of generality that
  \begin{align}\label{eq:d_0}
    d_0<\tilde\delta-\delta_1,
  \end{align}
  letting $T_i$ grow if necessary.

  In the opposite direction, we can find an upper bound for
  the distortion time associated with a given distance to
  $\cS_i$ reversing the inequality in~\eqref{eq:distS1} as
  follows. Letting $\gamma$ denote a smooth curve
  $\gamma:[0,1]\to M$ such that $\gamma(0)\in\cS_i$ and
  $\gamma(1)=x$ for any given fixed $x$ near $\cS_i$, we get
  \begin{align}
    \dist(f(x),f(\cS_i))
    &=
    \inf_\gamma \int_0^1
    \|Df(\gamma(t))\dot\gamma(t)\|\,dt
    \ge
    \frac1B\inf_\gamma
    \int_0^1 \dist(\gamma(t),\cS_i)^{\alpha_i}
    \|\dot\gamma(t)\|\,dt \nonumber
    \\
%     &=
%     \frac1B\inf_\gamma
%     \int_0^1\dist(\gamma(t),\cS)^\beta
%     \frac{d}{dt}\dist(\gamma(t),\cS)\frac{dt}
%     {|\cos\measuredangle(t)|} \nonumber
%     \\
    &\ge
    \frac1B \inf_\gamma
    \int_0^1 \dist(\gamma(t),\cS_i)^{\alpha_i}
    \left|\frac{d}{dt}\dist(\gamma(t),\cS)\right| \, dt
    \nonumber
    \\
    &=
    \frac1B \inf_\gamma
    \int_0^1\left|
      \frac{d}{dt}
      \frac{\dist(\gamma(t),\cS_i)^{1+\alpha_i}}{1+\alpha_i}
    \right| \, dt
    \nonumber
    \\
    &=
    \frac1{B(1+\alpha_i)} \inf_\gamma
      \var_{t\in[0,1]}\dist(\gamma(t),\cS)^{1+\alpha_i}
    \ge\label{eq:distS2}
    \frac{\dist(x,\cS)^{1+\alpha_i}}{B(1+\alpha_i)},
  \end{align}
  where we use, beside the non-flat condition (S1), the
  relation
  \begin{align*}
    \left|\frac{d}{dt}\dist(\gamma(t),\cS_i)\right|=
    \|\pi_t(\dot\gamma(t))\|=\|\dot\gamma(t)\|\cdot |\cos
    \measuredangle(\dot\gamma(t),N_t)|\le \|\dot\gamma(t)\|
  \end{align*}
  and write $N_t$ for the normal direction to the level
  submanifold
  \begin{align*}
    S_t=\{z\in M:\dist(z,\cS_i)=\dist(\gamma(t),\cS_i)\}
    \quad\text{at}\quad \gamma(t);
  \end{align*}
  and $\pi_t$ for the orthogonal projection from
  $T_{\gamma(t)}M$ to $N_t$. We also use the well known
  relation $ \var_{[0,1]} \vfi = \int_0^1 |D\vfi(t)|
  \,dt$ for the total variation of a differentiable function
  $\vfi:[0,1]\to\RR$.

  From the inequality \eqref{eq:distS2} we see that if
  $d_j\ge\dist(x,\cS)\ge d_{j+1}$ and $f(x)\in B_i$ for some
  hyperbolic neighborhood with distortion time $h$, then
  \begin{align*}
    \delta_2 e^{-\xi_0 c_0 h}
    &\ge
    \dist(f(x),f(\cS_i))
    \ge
    \frac{C_2^{1+\alpha_i}}{B(1+\alpha_i)}
    e^{ - \xi_0 c_0 (T_i+j+1)}
    =
    \frac{C_2^{1+\beta_i+\alpha_i-\beta_i}}{B(1+\alpha_i)}
    e^{ - \xi_0 c_0 (T_i+j+1)}
    \\
    &=
    \delta_2^{\alpha_i-\beta_i}
    \frac{(1+\beta_i)^{\alpha_i-\beta_i}}
    {B^{1+\alpha_i-\beta_i}(1+\alpha_i)}
    e^{ - \xi_0 c_0 (T_i+j+1)}
    \ge
    C_3 \delta_2^{\alpha_i-\beta_i} e^{- \xi_0 c_0 (T_i+j+1)}
    \quad \text{or}
    \\
    h
    &\le
    (T_i+j+1)
    - \frac{\log C_3}{\xi_0 c_0}
    +\frac{1-(\alpha_i-\beta_i)}{\xi_0 c_0}\log\delta_2
    \le
    (1+\frac\zeta2)(T_i+j+1)
  \end{align*}
  as long as $T_i$ is big enough, depending on $f$ and
  $\xi_0,c_0,\delta_2$. We remark that we have used the
  condition $0<\alpha_i-\beta_i<1$ in the inequalities
  above.  For future reference we write this inequalities
  (for a big enough $T_i$) in the convenient format
  \begin{align}
    \label{eq:distS3}
   d_j\ge\dist(x,\cS_i)\ge d_{j+1} , j\ge0
   \implies
   \frac{1+\beta_i}{\xi_0 c_0}\fD_{d_0}(x)\le h\le
   \frac{(1+\zeta/2)(1+\alpha_i)}{\xi_0 c_0}\fD_{d_0}(x).
  \end{align}
  We are now ready for the main arguments.

  \begin{claim}
    \label{claim:a}
    There are $(e^{\xi c_0},\delta)$-hyperbolic
    neighborhoods for each point in
    $B(\cS_i,\delta)\setminus\cS_i$, for suitable constants
    $\xi,\delta\in(0,1)$.
  \end{claim}

  Indeed, for each $y_0\in B(\cS_i,d_0)\setminus\cS_i$ there
  exists a unique integer $k$ such that $y_0\in
  B(\cS_i,d_k)\setminus B(\cS_i,d_{k+1})$. By the choice of
  $k$ we know that $y_1=f(y_0)$ has some distortion time
  $\frac{1+\beta_i}{1+\alpha_i}(T_i+k) \le h \le
  (1+\zeta/2)(T_i+k+1)$. Using the non-flat conditions on
  $\cS_i$ we can estimate the norm of the derivative and the
  volume distortion in a suitable neighborhood of $y_0$, and
  show that $y_0$ has a hyperbolic neighborhood with $h+1$
  as distortion time, with slightly weaker constants of
  expansion and distortion, as follows.

  Let $B_i$ be the hyperbolic neighborhood containing $y_1$
  with distortion time $h$ (as defined in
  \eqref{eq:halfball}).  The image of $B_i$ under $f^h$ is
  the $\delta_2$-ball around $z_h$ for some point $z_0\in
  f(\cS_i)$ and $y_{h+1}$ is inside this ball. Note that since
  $\dist(y_1,f(\cS_i))\ge
  d_{k+1}^{1+\alpha_i}/(B(1+\alpha_i))$
  from \eqref{eq:distS2}, then
  \begin{align}\label{eq:distB}
    \dist(y_{h+1},z_h)\ge
    \frac{d_{k+1}^{1+\alpha_i}}{B(1+\alpha_i)}
    e^{\xi_0 c_0 h}
    =
    \frac{C_2^{1+\alpha_i}}
    {B(1+\alpha_i)}e^{\xi_0 c_0 h
      (1- \frac{T+k+1}{h})}
    \ge  C_3 \delta_2^{\alpha_i-\beta_i}.
  \end{align}
  Thus if we set $2\delta_3=C_3
  \delta_2^{\alpha_i-\beta_i}$, then we can take a
  hyperbolic neighborhood of $y_1$ defined by $W:=(f^h\mid
  V_i)^{-1}(B(y_{n+1},\delta_3))$, where $V_i$ is the
  original neighborhood associated to $B_i$,
  see~\eqref{eq:halfball}. Observe that because $\diam W\le
  \delta_3 e^{-\xi_0 c_0 h}$
  \begin{align*}
    \dist(W,f(\cS_i))
    &\ge
    \frac{C_2^{1+\alpha_i}}{B(1+\alpha_i)}
    d_{k+1}^{1+\alpha_i}
    - \delta_3 e^{-\xi_0 c_0 h}
    \\
    &=
    \left(C_3\delta_2^{\alpha_i-\beta_i}d_{k+1}^{1+\alpha_i}
      e^{\xi_0 c_0 h} -\delta_3 \right)e^{-\xi_0 c_0 h}
    \ge
    \delta_3 e^{-\xi_0 c_0 h}
  \end{align*}
  by the definition of $\delta_3$.

Now we find a radius $\rho>0$ such that the image
$f(B(y_0,\rho))$ covers the hyperbolic neighborhood $W$ of
$y_1$. By the definition of $k$ and the non-flatness
condition (S1) we have that $f(B(y_0,\rho))\supset
B(y_1,\rho_1)$ for all small enough $\rho>0$, where
$\rho_1\ge \rho \|Df(y_0)^{-1}\|^{-1} \ge B^{-1}\rho
d_{k+1}^{-\alpha_i}$. Since we want $\rho_1\ge \delta_3
e^{-\xi_0 c_0 h}$ it is enough that
\begin{align*}
  \frac\rho{B} e^{-\frac{\alpha_i}{1+\alpha_i}\xi_0 c_0
    (T+k+1)}
  \ge
  \delta_3 e^{-\xi_0 c_0 h}
  \iff
  \rho\ge B\delta_3 e^{-\xi_0 c_0\frac{T+k+1}{1+\alpha_i}
    \big( (1+\alpha_i)\frac{h}{T+k+1} - \alpha_i\big)}.
\end{align*}
But using the relations obtained between $h$ and $T+k+1$ we get
\begin{align*}
  \rho\ge B\delta_3 e^{-\xi_0 c_0\frac{T+k+1}{1+\alpha_i}
    \big( (1+\alpha_i)(1+\zeta/2) - \alpha_i\big)}
  =
  B\delta_3 d_{k+1}^{1+(1+\alpha_i)\zeta/2}.
\end{align*}
% Note that by the non-flatness condition (S1) and the choice
% of $k$ we have $\|Df(w)^{-1}\|\le B d_{k+2}^{-\alpha_i}$
% for all $w\in M\setminus B(\cS_i,d_{k+2})$.
We note that $\gamma:=(1+\alpha_i)\zeta/2$ is positive due
to the non-flatness conditions.  If we take $w,\tilde w\in
B(y_0, d_{k+1}^{1+\gamma})$, then
\begin{align}\label{eq:V0-dentro}
  \dist(w,\cS_i)\ge\dist(y_0,\cS_i) - \dist(w,y_0)
  \ge d_{k+1}-d_{k+1}^{1+\gamma}
  = d_{k+1}\big(1-d_{k+1}^{\gamma}\big)
  >d_{k+2}
\end{align}
whenever $T_i$ is big enough in order that
$1-d_{k+1}^\gamma> d_{k+2}/d_{k+1}=e^{-\xi_0 c_0
  /(1+\alpha_i)}$.
% \begin{align*}
%   \left|\log \frac{\|Df(w)^{-1}\|}{\|Df(\tilde w)^{-1}\|}\right|
%     \le
%     B d_{k+1}^{-\bar\beta} \dist( w, \tilde w)
%     \le
%     \frac{B d_{k+1}}2
%     \quad\text{and}\quad
%     \|Df(w)^{-1}\|
%     \le
%     d_{k+1}^{-\bar\beta}
% \end{align*}
% by the condition (S2) on $f$  as long as we take $T$ big
% enough so that $\frac{e^{ B d_{k+1} }}{B}<1$.
Thus if $w_1=f(w) $ and $ \tilde w_1=f(\tilde w)$ are both
in $W$, then using the bound for the inverse of the
derivative provided also by the non-flatness condition (S1)
\begin{align}
  \dist(w,\tilde w)
  &\le
    B d_{k+2}^{-\alpha_i}\dist(w_1,\tilde w_1)
    \le
    B d_{k+2}^{-\alpha_i} e^{-\xi_0 c_0 h}\dist
    (w_{h+1},\tilde w_{h+1})\nonumber
    \\
  &=
  B\exp\left(-\xi_0 c_0 (h+1) \big(\frac{h}{h+1} -
    \frac{\alpha_i}{1+\alpha_i}\frac{T_i+k+2}{h+1}\big)\right)
  \dist(w_{h+1},\tilde w_{h+1})\nonumber
  \\
  &\le  \label{eq:xi-1}
  e^{-\xi_0 c_0 (1-\frac{\alpha_i}{1+\beta_i}) (h+1)}
  \dist(w_{h+1},\tilde w_{h+1})
\end{align}
for all big enough
$h$ since $h+1>\frac{1+\beta_i}{1+\alpha_i}(T_i+k)+1>
\frac{1+\beta_i}{1+\alpha_i}(1+o(T_i))(T_i+k+2)$ where
$o(T_i)$ is a small quantity which tends to zero when $T_i$
is taken arbitrarily large. We recall that condition (S1) on
$\alpha_i,\beta_i$ ensures that $1+\beta_i>\alpha_i$, so
that $0<\xi_1:=\xi_0\big(1-\frac{\alpha_i}{1+\beta_i}\big)$
and the contraction rate is bounded by $e^{-\xi_1 c_0}$. In
particular this shows that we can indeed take a neighborhood
with the radius $\rho$ we estimated above.%  since the
% derivative of the inverse is well defined and a
% contraction on the image of this neighborhood.  We can
% ensure that every distortion time $h$ satisfies this by
% taking $T_i$ sufficiently big.
% Moreover from~\eqref{eq:xi-1} we get
% \begin{align}\label{eq:xi-10}
%   B d_{k+1}^{-\alpha_i}\delta_3 e^{-\xi_0 c_0 h}
%   \le
%   \frac{ d_{k+1}^{1+\alpha_i}}2
%   \iff
%   B \delta_3 e^{-\xi_0 c_0 h}
%   \le
%   \frac{ d_{k+1}}2
%   =
%   e^{-\xi_0 c_0 h
%     \frac{T_i+k+1}{h(1+\alpha_i)}},
% \end{align}
% and, by the choice of $k$, the last expression on the right
% hand side of \eqref{eq:xi-10} is bounded from below by
% $e^{-\xi_0 c_0 h /(1+\beta_i)}$ which makes the inequalities
% in~\eqref{eq:xi-10} true as long as $h$ is big enough.

This means that every point $w_1$ in $W$ is the image of
some point $w\in B(y_0,d_{k+1}^{1+\gamma})$, and the
connected component $W_0$ containing $y_0$ of the pre-image
of $W$ under $f$ is fully contained in $M\setminus
B(\cS_i,d_{k+2})$, by the relation~\eqref{eq:V0-dentro}.
% Since we can take $\xi_0$ small enough so that $1-
% e^{-2b\xi_0 c_0/(1+\beta)} <1/2$, then the distance between
% any point in $V$ and $y$ is bounded by
%   \begin{align}\label{eq:innerdiam}
%     d_{k}-d_{k+2}=C e^{-b\xi_0 c_0 k/(1-\beta)} (1- e^{-2b\xi
%       c_0/(1-\beta)}) \le \frac12 \dist(y,\cS)
%   \end{align}
  Thus the neighborhood $W_0$ of $y_0$ satisfies items (1)
  and (2) of Proposition~\ref{p.contr} for $n=h+1$ iterates
  and for a $e^{-\xi_1 c_0}$ contraction rate.  Therefore item (3)
  of the same proposition also holds as a consequence with
  $\delta=d_0$, see Section~\ref{sec:hyperb-times}.

  In addition we can estimate
  \begin{align}
    S_{h+1}\fD_{d_0}(y_0)
    &=
    \fD_{d_0}(y_0) + S_h\fD_{d_0}(y_1)
    \le
    -\log d_{k+1} + \zeta h + o(\tilde\delta)\nonumber
    \\
    &\le
    \left( \xi_0 c_0 \frac{T_i+k+1}{(h+1)(1+\alpha_i)}
     +\zeta\frac{h}{h+1} + \frac{o(\tilde\delta)}{h+1}
   \right) (h+1)\nonumber
   \\
   &\le
   \left(\frac{\xi_0 c_0}{1+\alpha_i} + 2\zeta\right)(h+1)
   \le
   \frac{\xi_0 c_0}{1+\zeta} (h+1)  , \label{eq:SR}
   % \le\xi_0 c_0 (h+1)
  \end{align}
  where we have used item (3) of Proposition~\ref{p.contr}
  applied to the hyperbolic neighborhood $W$ of $y_1$ and
  the choice of $\zeta$ in~\eqref{eq:zeta}
  and~\eqref{eq:Daverage}. %Theorem~\ref{thm:abv0}.
  This bound is essential to obtain positive frequency of
  hyperbolic neighborhoods in the final stage of the
  proof. We use here the assumption B of the statement
  of Theorem~\ref{thm:Misiurew} with appropriately chosen
  constants.

  Moreover we also have the following estimate for the
  bounded distortion of volume, again using the non-flat
  conditions together with the above inequalities%  similarly
%   to \eqref{eq:xi-10},
  for $w,\tilde w\in W_0$
  \begin{align*}
    \log\frac{|\det Df(w)|}{|\det Df(\tilde w)|}
    &\le
    B \dist(W_0,\cS_i)^{-\alpha_i} \dist(w,\tilde w)
    \le
    e^{-\xi_1 c_0 (h+1)}\dist(w_{h+1},\tilde w_{h+1}).
  \end{align*}
  Therefore we can bound (using that $f(w),f(\tilde w)\in
  B_i$ and Proposition~\ref{p.contr})
  \begin{align*}
    \log\frac{|\det Df^{h+1}(w)|}{|\det Df^{h+1}(\tilde w)|}
    &\le
    \sum_{j=0}^{h} \log\frac{|\det Df(w_j)|}{|\det
      Df(\tilde w_j)|}
    \\
    &\le
    \sum_{j=0}^{h} e^{-\xi_1 c_0 j}
    \dist(w_{h+1},\tilde w_{h+1} )
    \le \frac{\delta_3}{1-e^{-\xi_1 c_0}} = C_1.
  \end{align*}
  Hence $W_0$ satisfies all the conditions  of
  Proposition \ref{p.contr} with appropriate constants.

  We stress that the value of $T_i$ (thus the value of
  $d_0$) depends only on $\alpha_i,\beta_i$ and
  $\xi_0,c_0,\delta_1$.

  This completes the proof of Claim~\ref{claim:a}.

  \begin{remark}
    \label{rmk:highexpnearsing}
    The minimal backward contraction in $B(\cS_i,d_0)$ near
    a singularity can be made arbitrarily big in a single
    iterate taking $T_i$ large enough or, which is the same,
    taking the neighborhood of the singularities small enough.
  \end{remark}

  \begin{remark}\label{rmk:exp-const}
    The constant of average expansion $\xi c_0$ is the same
    in all cases for each connected component of $\cS$, but
    in principle the radius $\delta_3=\delta_3(i)$ and the
    distance $d_0=d_0(i)$ to $\cS_i$ ensuring existence of
    distortion times depend on the connected
    component. However we assume that the values of $T_i$
    are big enough so that all the inequalities above are
    satisfied for all connected components $\cS_i$ and such
    that the corresponding neighborhoods of $\cS_i$ be
    contained in the $\tilde \delta$-neighborhood of $\cS$.
  \end{remark}

  In what follows we write $\delta$ (which is smaller than
  $\tilde\delta$) for the smallest value $d_0(i)$ over all
  connected components.

\subsubsection{Hyperbolic neighborhoods for almost every point}
\label{sec:hyperb-neighb-almost}

We use the statement of Claim~\ref{claim:a} from now on
and show that all points whose orbit does not fall into the
singular/critical set admit some distortion time with well
defined contraction rate and slow approximation rate.

\begin{claim}
 \label{claim:b}
 Lebesgue almost every point admits a $(e^{-\xi_2
   c_0},\delta)$-hyperbolic neighborhood for some iterate
 $\ell$ and the frequency of visits of the $\ell$-iterates
 to a $\delta$-neighborhood of $\cS$ is bounded by
 $\frac{\xi_0c_0}{1+\zeta}$.
\end{claim}

In a similar way to the one-dimensional multimodal case, we
use condition D on $\hat U$ to show that we can obtain a
minimal derivative when an orbit returns to any fixed
arbitrarily small neighborhood $U\subset\hat U$ of the
critical/singular set. From this intermediate result we
deduce that most points on $M\setminus U$ have some
distortion time.  We prove first an auxiliary result.

\begin{claim}
  \label{claim:b1}
  There exists a minimal average expansion rate $e^{- c/2}$
  either for the first return of $x_0$ to a small enough
  neighborhood $U$ of $\cS$, or for the first $N$
  iterates, where $N$ does not depend on the starting
  point $x_0\in U$.
\end{claim}

Let $U=B(\cS,\delta)$ be an open neighborhood of $\cS$
compatible with the choices of the $T_i$ as in the proofs of
the previous claims. So $U$ is the union of a number of
connected open sets, one for each connected component of
$\cS$.

The first $h$ iterates of $x_0$ correspond to a $(e^{-\xi
  c_0},\delta)$-distortion time.  From condition (C) there
are $K,c>0$ such that for $n\ge h$ satisfying $x_h,\dots,
x_{n-1}\in M\setminus U$ we have that $S_n\psi(x_0)\le -\xi
c_0 h + K - c (n-h)$.
 
We shrink the neighborhood $U$ so that the smaller
distortion time $\underbar h$ for points $x_0\in U\setminus\cS$
satisfies
\begin{align}\label{eq:cond-h0}
  \frac{\kappa}{\underbar h}
  <
  \min\left\{\frac{\xi c_0}4,\frac{c}8\right\}
  \quad\text{and}\quad
  \frac{\underbar h}{\underbar h+2 K/ c}\ge \frac12.
\end{align}
We write $[t]:=\max\{i\in\ZZ:i\le t\}$ for all $t\in\RR$ in
what follows and set $N=[2 K/ c]$.
\begin{description}
\item[Case (i) -- the orbit returns to $U$ in more than $N$
  iterates] In this case we have for $n > N$ that $
  S_n\psi(x_0) \le -\frac{ c}2 n$,
  since we can assume without loss that $c < \xi c_0$.
  Notice that we have the same conclusion if the orbit of
  $x_0$ \emph{never returns} to $U$.
\item[Case (ii) -- the return to $U$ occurs in less than
  $N$ iterates] We now use condition D to get, since the
  first return iterate $n$ satisfies $n\ge h\ge \underbar h$
  \begin{align*}
    S_n\psi(x_0)\le -\xi c_0 h +\kappa
    =
    n (-\xi c_0\frac{h}n + \frac\kappa{n})
    \le
    n (-\frac{\xi c_0}2 + \frac{\xi c_0}4)
    = -\frac{\xi  c_0}4 n,
  \end{align*}
  since $\underbar h$ was chosen as in~\eqref{eq:cond-h0}.
\end{description}
This completes the proof of Claim~\ref{claim:b1}, setting
$\min\{ c/2, \xi c_0/4\}$ as the average expansion rate.

% \begin{description}
% \item[blabla]
%   Divide the iterates into segments of orbit in distortion
%   times starting at the returns $n_i$ to $\hat U$:
%   \begin{align*}
%     \underbrace{x_0,\dots,x_{h_0}}_\text{hyp. time}
%     \underbrace{,\dots,}_\text{in $M\setminus \hat U$}
%     \underbrace{x_{n_1},\dots,x_{n_1+h_1}}_\text{hyp. time}
%     \underbrace{,\dots,}_\text{in $M\setminus \hat U$}
%     \underbrace{x_{n_2},\dots,
%     x_{n_2+h_2}}_\text{hyp. time}, \dots
%   \end{align*}

%   \begin{remark}
%   \label{rmk:smallernearsing}
%   The constant $\hat K$ depends only on the components of
%   $\hat U$ near the critical set $\cC:=\{x\in\cS: \det
%   Df(x)=0\}$ by the non-flat conditions on $\cS$, since near
%   singularities the derivative is strongly expanding and the
%   manifold is compact.
% \end{remark}
% Hence we may shrink $\hat U$ near $\cS\setminus\cC$ without
% affecting $\hat K$, but making the minimum contraction
% $\eta$ through singular points stronger than both $\hat
% K=\hat K(\hat U)$ and $K=K(U)$ (both given by condition (C))
% in what follows.

Now we prove Claim~\ref{claim:b}.

Again fix an arbitrary $x_0\in M\setminus U$ such that the
orbit of $x_0$ never falls into $\cS$. We divide the
argument in two cases.

\begin{description}
\item[Case (iii) -- the orbit takes more than $N$ iterates
  to enter $U$] Arguing as in the proof of Case (i) above we
  get for all $n > N$ satisfying $x_0,\dots, x_{n-1}\in
  M\setminus U$ that $ S_n(x_0)\le K - cn \le -\frac{c}2 n.$
  This implies that there exists some distortion time $h\le
  n$ for $x_0$ by the proof of Theorem~\ref{thm:abv0}. We
  note that in this case $S_n\fD_{\delta}(x_0)=0$.

\item[Case (iv) -- the orbit enters $U$ in at most $N$
  iterates] Let $j\le N$ be the first entrance time of the
  orbit of $x_0$ in $U$, $h$ the distortion time associated
  to $x_j$ and use Claim~\ref{claim:b1} to get
\begin{align*}
  S_{j+h}\psi(x_0)
  &\le
  \kappa + S_{h}\psi(x_{j})
  \le
  \kappa- \frac{c}2 h
  =
  (j+h)\Big( \frac{\kappa}{j+h}
  - \frac{c}2\frac{h}{j+h} \Big)
    \\
  &\le
  (j+h)\Big( \frac{c}8 - \frac{c}4\Big)=-\frac{c}8 (j+h)
\end{align*}
since $h\ge\underbar h$ by the choice of $U$ in
\eqref{eq:cond-h0}.

We also obtain that
$S_{j+h}\fD_{\delta}(x_0)=S_h\fD_\delta(x_j)\le
\frac{\xi_0c_0}{1+\zeta} h\le \frac{\xi_0c_0}{1+\zeta} (j+h)$.

Moreover the distortion of $f^{j+h}$ on the connected
component of $(f^j)^{-1}(V_{x_j})$ containing $x_0$ is
bounded from above by a constant dependent on $N$ only
(since we are dealing with a local diffeomorphism away from
a given fixed neighborhood of $\cS$), where $V_{x_j}$ is the
hyperbolic neighborhood of $x_j$ given by
Claim~\ref{claim:a}.
\end{description}

This completes the proof of Claim~\ref{claim:b} if we set
$\xi_2$ such that $\xi_2 c_0 = c/8$.

%%%%%%%%%%%%%%%%%%%%%%%%%%%%%%%%%%%%%%%%%%%%%%%%%%%%%%%%%%%

\subsubsection{Positive frequency of hyperbolic
  neighborhoods almost everywhere}
\label{sec:positive-frequency-h}

Here we finish the proof of Theorem~\ref{thm:dense-hyp-times}.

\begin{claim}
  \label{claim:c}
  The frequency of $(e^{-\xi_2 c_0},\delta)$-hyperbolic
  neighborhoods is positive and bounded away from zero
  Lebesgue almost everywhere.
\end{claim}

We now define an auxiliary induced map $F:\tilde M\to
M$, so that the $F$-iterates correspond to iterates of $f$
at distortion times, as follows. For $x_0\in\tilde M$ we
have two cases
  \begin{description}
  \item[Expansion without shadowing] the $\ell\ge N$
    iterates of $x$ belong to $M\setminus U$. In this case
    there exists some $(e^{-\xi_2 c_0},\delta)$-distortion
    time $\ell\le N$ and we define $F(x)=f^\ell(x)$ and
    $\tau(x)=\ell=q(x)$.
  \item[Expansion with shadowing] let $0\le q<N$ be the
    least non-negative integer such that $x_q\in U$ and $p$
    be the $(e^{-\xi_0 c_0},\delta)$-distortion time
    associated to $x_q$ from Claim~\ref{claim:a}. We define
    $F(x)=f^{p+q}(x)=x_{p+q}$ and $\tau(x)=p+q$ and $q(x)=q$
    in this case.
  \end{description}
  The images of $F$ always belong to $M\setminus
  B(\cS,\delta)$ by the choice of $\delta$
  from~\eqref{eq:d_0} following
  Remark~\ref{rmk:preball-away}.  In addition, for any given
  $x_0\in\tilde M$ the map $F(x)$ is defined and the iterate
  $\tau(x)$ has all the properties of a $(e^{-\xi_2
    c_0},\delta)$-distortion time for $x_0$ according to
  Claim~\ref{claim:b}.

  Now we observe the statement of Claim~\eqref{claim:c}
  is a consequence of the following property: there exists
  $\theta>0$ such that for every $x\in\tilde M$
  \begin{align}\label{eq:positive-freq}
    \limsup_{n\to+\infty}\frac1n\#\{ 0\le j < n :
    f^j(x)\in\cO^+_F(x)\} \ge\theta,
  \end{align}
  where $\cO^+_F(x)=\{F^i(x), i\ge0\}$ is the positive orbit
  of the induced map $F$. Moreover it is easy to see that
  for each $x\in M$ and each $n\in\NN$
  \begin{align}\label{eq:updown}
    \#\{ 0\le j < n : f^j(x)\in\cO^+_F(x)\}
    =
    \sup\big\{
    k\ge0: S_{k+1}^F\tau(x)=\sum_{i=0}^k
    \tau\big(F^i(x)\big) < n
    \big\}.
  \end{align}
  We remark that from \eqref{eq:updown} to obtain
  \eqref{eq:positive-freq} it is enough to show that
  $S_{k+1}^F \tau(x)< \frac1{\theta}k$, at least for every big
  enough $k$. Indeed
  $\big\{
    k\ge0: S_{k+1}^F\tau(x) < n
    \big\}
    \supset
    \big\{
    k\ge0: \frac1\theta k < n
    \big\}$
   so
   $ \sup\big\{
    k\ge0: S_{k+1}^F\tau(x) < n
    \big\} \ge \theta n.$

  The bounds~\eqref{eq:distS3} and~\eqref{eq:SR} together
  with the definition of $F$ ensure that we are in the
  conditions of the following result.

  \begin{lemma}
    \label{le:positive-freq}
    Assume that we have an induced map $F=f^\tau$ for some
    $\tau:K\to\NN$ defined on a positive invariant subset
    $K$ such that for every $x\in K$:
    \begin{enumerate}
    \item $\tau(x)=q(x)+p(f^{q(x)}(x))$ for well defined
      integer functions $q$ on $K$ and $p$ on $f^{q(x)}(x)$
      for all $x\in K$;
    \item there exists $N\in\NN$ such that $q\le N$;
    \item there exists $0<d<1$ and $C,\rho>0$ such that
      $0<\rho C<1$ and for all $x\in K$
      \begin{enumerate}
      \item the iterates $x,f(x),\dots,f^{q(x)-1}(x)$ are
        outside $B_{d}(\cS)$;
      \item $p(f^{q(x)}(x)) \le C \fD_d(f^{q(x)}(x))$;
      \item $S_p^f \fD_d(f^{q(x)}(x))\le\rho p$.
      \end{enumerate}
    \end{enumerate}
    Then $S_k^F\tau(x)\le\frac{N}{1-\rho C} k$ for each
    $x\in K$ and every $k\ge1$.
  \end{lemma}

  In addition, we observe that since $f$ is a regular map
  (that is $f_*\Leb\ll\Leb$) then we can further assume that
  the full Lebesgue measure set $\tilde M$ is $f$-invariant,
  since $\cap_{i\ge0}f^{-i}(\tilde M)$ also has full Lebesgue
  measure. Hence we can apply Lemma~\ref{le:positive-freq}
  with $K=\tilde M$, $C=\frac{1+\zeta/2}{\xi_0 c_0}$
  and $\rho=\frac{\xi_0 c_0}{1+\zeta}$ to obtain
  $\theta\ge\frac{\zeta}{2(1+\zeta)\ell}$.

  \begin{remark}\label{rmk:log-dist-integrable}
    According to item (3b) of the statement of
    Lemma~\ref{le:positive-freq} above, if we have an
    absolutely continuous invariant probability measure
    $\nu$ for $f$, then $\nu(\fD_d)\le\rho<\infty$ and so
    $\fD_d$ is $\nu$-integrable, as claimed in
    Theorem~\ref{thm:Misiurew}.
  \end{remark}

  This completes the proof of
  Theorem~\ref{thm:dense-hyp-times} except for the proof of
  Lemma~\ref{le:positive-freq}.

\begin{proof}[Proof of Lemma~\ref{le:positive-freq}]
  Using the definition of $F$ and the assumptions on $\tau$,
  for every given $k\ge0$ and $x\in K$ we can associate a
  sequence $q_0,p_0,q_1,p_1,\dots,q_k,p_k$ such that for
  each $i=0,\dots,k$ we have $ q_i=q(F^i(x))
  \quad\text{and}\quad q_i+p_i=\tau(F^i(x)).$ This
  together with the assumptions of item (3) in the statement
  of the lemma allows us to estimate
  \begin{align}
    S_{k}^F\tau(x) &= \sum_{i=0}^k ( q_i + p_i) \le
    \sum_{i=0}^{k-1}\Big( N +
    C\fD_{d_0}\big(f^{q_i}(F^i(x))\big) \Big) \nonumber
    \\
    &\le k N + C\sum_{i=0}^{k-1}
    S^f_{q_i+p_i} \fD_{d_0} (F^i(x)) \label{eq:noreturn}
    \\
    &\le \label{eq:contraction1} k N
    + C\rho
    \sum_{i=0}^{k-1}(q_i+p_i) = k N +
    C\rho S_{k}^F\tau(x),
  \end{align}
  where in~\eqref{eq:noreturn} we have used that
  $\fD_{d_0}\ge0$ and that this function equals zero at each
  of the $q_i$ iterates before each visit to $B(\cS,d)$.
  The contraction in~\eqref{eq:contraction1} implies $
  S_k^F\tau(x) \le \frac{N}{1-\rho C}k $ as stated.
\end{proof}

%%%%%%%%%%%%%%%%%%%%%%%%%%%%%%%%%%%%%%%%%%%%%%%%%%%%%%%%%%%%%

%%%%%%%%%%%%%%%%%%%%%%%%%%%%%%%%%%%%%%%%%%%%%%%%%%%%%%%%%

\section{Periodic attractor with full basin of attraction}
\label{sec:full-basin-attract}

Here we prove item (3) of Theorem~\ref{mthm:codimension}. We
show that a perturbation of $X$, like the one depicted in
the right side of Figure~\ref{fig:skewperturb}, for a flow
$Y\in\cP\cap\cN$ is such that $U$ is as a trapping region
which coincides Lebesgue modulo zero with the basin of a
periodic attracting orbit (a sink) of $Y$. This is a
consequence of the smoothness of the first return map to
$\Sigma_1$ after quotienting out the stable leaves, which is
the reason why we assume the flow is of class at least $C^2$
and restrict the vector field to a submanifold $\cP\cap\cN$
of all possible vector fields nearby $X$, together with the
robustness of property C in the statement of
Theorem~\ref{thm:Misiurew}.

Indeed, let $g=g_Y$ be the action on the stable leaves of
the first return map $R_Y$ of the flow of $Y\in\cP\cap\cN$
to $\Sigma_2$. Recall that we constructed $X$ having on
$\Lambda_\Sigma=\cap_{n\in\ZZ^+}g^n(\Sigma_2)$ a partially
hyperbolic splitting so that the stable foliation does
persist under perturbations.  The projection along the
leaves of this foliation in $\Sigma_2$ is absolutely
continuous with H\"older Jacobian, as a consequence of the
strong domination obtained in
Section~\ref{sec:dominat-splitt}.

% If the foliation has codimension one in
% $\Sigma_2$, then the projection is in fact H\"older-$C^1$,
% as explained in \cite[Section 5]{APPV}.

% However the stable foliation now no longer has codimension
% one. But still this foliation is absolutely continuous, as
% explored by Anosov~\cite{An67} when proving ergodicity of
% the geodesic flow for manifolds with negative sectional
% curvature with dimension higher than two. In fact, the
% Jacobian of the holonomy maps along the leaves of any
% uniformly expanding or contracting foliation is H\"older
% continuous with respect to the Riemannian metric on the
% manifold, see e.g. \cite{Vi97b,barreira-pesin2,PSW97}. More
% precisely, we note that the time 1 map $Y^1$ of $Y$ is also
% a partially hyperbolic diffeomorphism with a strong stable
% direction $E^{s}$. We recall that we stay in the set of
% $C^2$ vector fields. By Theorem A in \cite{PSW97}, applied
% to $Y^1$, the (strong) stable holonomy map is H\"older
% continuous. Note that the transversal section $\Sigma_2$ for
% $X$ is also transversal for $Y$, as we take $Y$ sufficiently
% close to $X$.  Therefore, there exists a neighborhood $V_2$
% of $\Sigma_2$ and a $C^2$-projection along the field
% trajectories $\eta: V_2 \to \Sigma_2$. So the stable
% foliation for the Poincar\'e map is H\"older, as it can be
% obtained by composing $\eta$ with the leaves in the strong
% stable foliation for $Y^1$ contained in $V$.  This is what
% we need to show that the map $g_Y$ is such that $\log|\det
% Dg_Y|$ is an H\"older function on $N$.

The map $g$ sends a $\delta$-neighborhood of $\cS$ into the
local basin of attraction of a periodic sink $p$. We denote
by $B$ the stable set of the orbit of $p$. On $N\setminus B$
the map $g$ is uniformly expanding, since $B\supset
B(\cS,\delta)$ and condition C on the statement of
Theorem~\ref{thm:Misiurew} is persistent under small
perturbations.

\begin{lemma}
  \label{le:fullbasinsink}
  Lebesgue almost every point of $N$ belongs to the basin of
  the periodic sink.
\end{lemma}

\begin{proof}
  Arguing by contradiction, assume that the basin $B$ of the
  sink is such that $E:=N\setminus B$ has positive
  volume. Since $B$ contains a neighborhood of $\sigma_1$
  and a neighborhood of the sink, then $E$ is an invariant
  subset satisfying condition C of the statement of
  Theorem~\ref{thm:Misiurew}. Hence $f_1\mid E$ is uniformly
  expanding: there exists $N\ge1$ and $\lambda\in(0,1)$ such
  that $\|(Dg_Y^N)^{-1}\|\le\lambda$. Therefore, since $g_Y$
  is $\log$-H\"older and expanding, and $E$ is closed,
  invariant and has positive volume, we can apply the
  arguments in \cite{AAPP} to show that there exists a ball
  $U$ of radius $r>0$ fully contained in $E$.

  We claim that for $g=g_Y^N$ there exists $\rho>1$ such
  that $g^k(U)$ contains a ball of radius $\rho^k r$ for all
  $k\ge1$, which yields a contradiction, since the ambient
  manifold is compact and $E$ is by assumption a proper
  subset.

  To prove the claim, recall that $g$ is a local
  diffeomorphism on a neighborhood of $E$, since $E$ is far
  from the singularities of the stable foliation. We assume
  that $B(x_0,s_0)$ is a ball centered at $x_0$ with radius
  $s_0$ contained in $E$ and consider $g(B(x_0,s_0))$.

  Let us take $y_1$ in the boundary of $g(B(x_0,s_0))$ and a
  smooth curve $\gamma_1:[0,1]\to N$ such that
  $\gamma_1(0)=x_1:=g(x_0)$ and $\gamma_1(1)=y_1$.  Let
  $\gamma_0$ be a lift of $\gamma_1$ under $g$, that is,
  $\gamma_1=g\circ\gamma_0$ such that $\gamma_0(0)=x_0$.  We
  then define $s:=\sup\{ t\in[0,1] : \gamma_1([0,t])\subset
  g(B(x_0,s_0)) \}.$ Clearly $s>0$ and by its definition and
  the expansion properties of $g$ in $E$ we get
  \begin{align*}
    \lambda \times (\text{length of  } \gamma_1([0,s]))
    \ge \text{length of  } \gamma_0([0,s]) \ge
    \dist(\gamma_0(0),\gamma_0 (s)).
  \end{align*}
  However, $\gamma_1(s)$ is at the boundary of
  $g(B(x_0,s_0))$ so that $\gamma_0(s)$ is also at the
  boundary of $B(x_0,s_0)$, because $g$ is a local
  diffeomorphism. Thus we get
  \begin{align*}
    \dist(y_1,x_1)\ge \text{length of } \gamma_1([0,s]) \ge
    \frac1\lambda\times \dist(\gamma_0(0),\gamma_0 (s)) \ge
    \frac1\lambda s_0
  \end{align*}
  and the claim is proved with $\rho=\lambda^{-1}$.
\end{proof}

%%%%%%%%%%%%%%%%%%%%%%%%%%%%%%%%%%%%%%%%%%%%%%%%%%%%%

\appendix

\section{Non-uniform expansion and existence of hyperbolic
times}
\label{sec:hyperb-neighb-constr}

Here we prove Theorem~\ref{thm:abv0}. The proof is very
similar to \cite[Lemma 5.4]{ABV00} but our definition of
hyperbolic times/hyperbolic neighborhoods
is slightly different, in a crucial way, from the definition
on \cite{ABV00}, and we include a proof for completeness.

We start with the following extremely useful technical
result will be the key for several arguments.

\begin{lemma}
  \label{le:Pliss}
  Let $H\ge c_2 > c_1 >0$ and
  $\zeta={(c_2-c_1)}/{(H-c_1)}$. Given real numbers
  $a_1,\ldots,a_N$ satisfying
 $$
 \sum_{j=1}^N a_j \ge c_2 N \quad\text{and}\quad a_j\le H
 \;\;\mbox{for all}\;\; 1\le j\le N,
 $$
 there are $l>\zeta N$ and $1<n_1<\ldots<n_l\le N$ such
 that $$ \sum_{j=n+1}^{n_i} a_j \ge c_1\cdot(n_i-n)
 \;\;\mbox{for each}\;\; 0\le n < n_i, \; i=1,\ldots,l. $$
\end{lemma}

\begin{proof}
  See \cite[Lemma 11.3]{Man87}.
\end{proof}

\begin{proof}[Proof of Theorem~\ref{thm:abv0}]
The proof uses Lemma~\ref{le:Pliss} twice,
first for the sequence $a_j=-\psi(x_{j-1})$
(properly cut off so that it becomes bounded from above),
and then for $a_j=\fD_\delta(x_j)$ for an adequate $\delta>0$.

Let $H\subset M\setminus\cS$ be such that conditions
\eqref{liminf1} and \eqref{eq:slow-recurrence} hold for all
$x\in H$ and let $0<\xi<1$ and $\zeta>0$ be given, and take
$x=x_0\in H$ and $\gamma_0:=(2+\xi)/3\in(\xi,1)$,
$\gamma_2:=(1-\xi)/3$ and
$\gamma_3=\gamma_1-\gamma_2=(1+2\xi)/3$.  Then for every
large $N$ we have $S_N\psi(x) \le - \gamma_0 c N.$ Moreover
since $f(\cS)\cap\cS=\emptyset$ we can assume that
$\zeta<\inf\{\dist(x,y): x\in\cS, y\in f(\cS)\}$.

For any fixed $\rho>\beta$, by non-degeneracy condition (S1),
we can find a neighbourhood $V$ of $\cS$ such that
\begin{align}
\label{eq:bound2}
|\psi(z)|
\le \rho \fD(z)\quad\text{for every}\quad
x\in V.
\end{align}
Setting $\varepsilon_1>0$ such that $\rho\varepsilon_1
\le \gamma_1$, we can use the slow approximation condition to
find $r_1>0$ so small that
\begin{align}
  \label{eq:bound1}
S_N\fD_{r_1}(x) \le \varepsilon_1 N.
\end{align}
We may assume without loss that $V=B(\cS,r_1)$ in what
follows.  Now we fix $H_1 \ge \max\{ c,\rho |\log r_1|,
\sup_{M\setminus V} |\psi| \}$ and define the set $E=\{z\in
M: \psi(z)<-H_1\}$ and the sequence
$a_j=-(\psi\chi_{M\setminus E})(x_{j-1})$.

\emph{We remark that there is a shift between the index of
  $a_j$ and that of $x_{j-1}$ in the above definition. }

We note that by construction $a_j\le H_1$ and that $x_j\in
E$ implies $x_j\in V$, $\fD(x_j)<-\log r_1$, because $\rho
|\log r_1| \le H_1 < -\psi(x_{j}) < \rho \fD(x_j).$

This means that $\fD_{r_1}(x_j)=\fD(x_j)<|\log r_1|$
whenever $x_j\in E$. From \eqref{eq:bound2}
and~\eqref{eq:bound1} we get that
$
  -S_N(\psi\chi_E)(x) \le \rho S_N(\fD\chi_E)(x)\le
  \rho\epsilon_1 N \le \gamma_1 c N.$

Therefore $\sum_{j=1}^{N} a_j = -S_N\psi(x) -
S_N(\psi\chi_E)(x) \ge (\gamma_0-\gamma_1)c N = \gamma_3 c
N$. Hence we can apply Lemma~\ref{le:Pliss} with
$c_2=\gamma_3 c$, $c_1=\xi c$, $H=H_1$, obtaining
$\theta_1=\gamma_2 c/(H-\xi c)\in(0,1)$ and $l_1\ge \theta_1
N$ times $1 \le p_1< \cdots <p_{l_1}\le N$ such that
  \begin{align}
    \label{eq:conclusion1}
    \sum_{j=n+1}^{p_i} \psi(x_{j-1}) \le -\sum_{j=n+1}^{p_i}
    a_j \le -\xi c  (p_i-n)
  \end{align}
for every $0 \le n < p_i$ and $1\le i\le l_1$.

Let now $\varepsilon_2>0$ be small enough so that
$\varepsilon_2 < \min\{\zeta, b c \theta_1\}$, and let
$r_2>0$ be such that $-S_N\fD_{r_2}(x_1) \ge -\epsilon_2 N$
from the slow recurrence condition (we note that
$\dist(x_0,\cS)>\zeta>r_2$).  Taking $c_1= b c$,
$c_2=-\varepsilon_2$, $A=0$, and
$
\theta_2=\frac{c_2-c_1}{A-c_1}=1-\frac{\varepsilon_2}{bc}
$
we can apply again Lemma~\ref{le:Pliss} to
$a_j=-\fD_{r_2}(x_j)$.

\emph{We remark that now there is no shift
between the index of $a_j$ and $x_j$ in the previous
definition. }

In this way we obtain $l_2\ge\theta_2 N$ times $1\le q_1 <
\cdots < q_{l_2}\le N$ such that
\begin{align}
  \label{eq:conclusion2}
  \sum_{j=n+1}^{q_i} \fD_{r_2}(x_j)\le \epsilon_2
   (q_i - n) \le \zeta (q_i-n)
\end{align}
for every $0 \le n < q_i$ and $1\le i \le l_2$.

Finally since our choice of $\theta_2$ ensures that
$\theta=\theta_1+\theta_2-1>0$, then there must be
$l=(l_1+l_2-N)\ge \theta N$ and $1\le n_1 <\ldots <n_l \le
N$ such that \eqref{eq:conclusion1} and
\eqref{eq:conclusion2} simultaneously hold.

This exactly means that we have condition~\eqref{d.ht} with
$\delta=r_2$, because $\fD_{r_2}(x_0)=0$ by the choice
of $r_2$ above, and $\xi c$ as the logarithm of the
contraction rate, as in the statement of
Theorem~\ref{thm:abv0}.
\end{proof}

\section{Solenoid by isotopy}
\label{sec:isotopy}

We recall that $\sS^1=\{z\in\CC: |z|=1\}$, $\TT=(\sS^1)^k$,
$B^k:=\{ x=(x_1,\dots,x_k)\in\RR^n: \sum_{i=1}^k x_i^2=1\}$
for all $k\ge1$, and $\cT^k=\TT^k\times\DD$, where
$\DD=\{z\in\CC: |z|<1\}$.

Here we prove the results needed in
Section~\ref{sec:unpert-singul-flow} ensuring the existence
of a smooth family of embeddings of $\cT^k$ into $B^{k+2}$,
for all $k\ge1$, which deforms a tubular neighborhood in
$B^{k+2}$ of the usual embedding of $\TT^k$ into
$B^{k+1}\simeq B^{k+1}\times\{0\}\subset B^{k+2}$, into the
embedding of $\cT$ into the image of the Smale solenoid map.

More precisely, consider the identity map $i$ on $\cT$ and
the solenoid map
\begin{align*}
  s:\cT\to\cT, \quad
  (\Theta=(\theta_1,\dots,\theta_k),z)\mapsto
  (\Theta^2=(\theta_1^2,\dots,\theta_k^2), A_\Theta(z))
%\lambda z + \rho  e^{i\sum_{l=1}^k 2^{-l}\theta_l})
\end{align*}
where $A_\Theta$ is a contraction with contraction rate
bounded by $0<\lambda<1$, and the map in the coordinate
$\Theta$ is the expanding torus endomorphism $f_1$ defined
in Section~\ref{sec:an-example-nue}, but restricted to
$\TT^k\subset\CC^k$.

\begin{proposition}\label{pr:embedding}
  There exists an embedding $e:\cT\to B^{k+2}$ such that the
  projections $\pi_\DD:\cT\to\TT^k$ on the first coordinate
  and $\tilde\pi_\DD:e(\cT)\to e(\TT^k\times\{0\})$ along
  the leaves of the foliation $\cF^s:=\{ e(\Theta\times\DD)
  \}_{\Theta\in\TT^k}$ of $e(\cT)$ define the solenoid map
  $S$ by the commutative diagram:
  \begin{center}
    $\begin{CD} \cT @>>{s}> \cT
      \\
      @VV{e}V @V{e}VV
      \\
      e(\cT) @>{S}>> e(\cT)
    \end{CD}.$
  \end{center}
  Moreover there exists a smooth family $e_t:\cT^k\to B^{k+2}$
  of embeddings for all $t\in[0,1]$ such that $e_0=e\circ i$
  and $e_1=e\circ s$.
\end{proposition}

We prove this statement in the following steps.

\subsubsection{An embedding of $\TT^k$ in $B^{k+2}$}
\label{sec:embedd-ttk-bk+1}

We argue by induction on $k\ge1$. We know how to embed
$\TT^1$ in $B^3$. Let us denote by $e^1$ this embedding
and fix a small number $d>0$ and $\lambda\in(0,1/2)$.

We assume that we have an embedding $e^{l}$ of $\TT^{l}$ on
$B^{l+2}$ for all $l=1,\dots,k-1$, in such a way that the
image of $e^{l+1}$ is in a tubular neighborhood of the image
of $e^l$ inside $B^{l+2}$ with size $\le d\lambda^{l+1}$, for
$1\le l < k-1$.

For each $w\in e^{k-1}(\TT^{k-1})\subset \RR^{k+1}\simeq
\RR^{k+1}\times 0 \subset \RR^{k+2}$ let
$N_w=(T_w\sS_{k})^\perp$ be the normal space to
$e^{k-1}(\TT^{k-1})$ at $w$ in $\RR^{k+2}$ (where we take in
$\RR^{k+2}$ the usual Euclidean inner product). This is a
$3$-dimensional space.

We know we can embed $\TT^1$ into $B^3$ through
$e^1$. So by a simple rescaling we can assume that $e^1_w$
embeds $\TT^{1}$ into a small neighborhood of $w$ in
$w+N_w$. To keep the inclusion in the tubular
neighborhood of the image of $e^{k-1}$, we take this
neighborhood around $w$ to have radius $\lambda^{k+1}d$.

Hence letting $\hat w\in\TT^{k-1}$ be the unique element
such that $e^{k-1}(\hat w)=w$ and considering the map
$e^k:\TT^{k-1}\times\TT^1\to B^{k+2}$ given by $(\hat
w,\theta)\mapsto e^1_w(\theta)$ we easily see that
\begin{itemize}
\item $D_1 e^k (\hat w,\theta)  (\RR^{k-1}) =
  De^{k-1}(\hat w)(\RR^{k-1})$ is the tangent space of
  $e^{k-1}(\TT^{k-1})$ at $w$;
\item $D_2 e^k (\hat w,\theta)  (\RR)$ is a subspace of $N_w$.
\end{itemize}
Therefore the tangent map $De^k$ to $e^k$ always has maximal
rank and clearly is injective with a compact domain, thus
$e^k$ is an embedding.  This completes the induction
argument and proves the existence of an embedding
$e^k:\TT^k\to B^{k+2}$.

\begin{remark}\label{rmk:embedd-B2}
  This argument is also true if we start with an embedding
  of $\TT^1$ into $B^2$ so that we obtain an embedding of
  $\TT^k$ into $B^{k+1}$ for all $k\ge1$. But we need one
  extra dimension to deal with the solid torus in what
  follows.
\end{remark}

\emph{It is crucial to observe that the entire inductive
  construction just presented is built over nested tubular
  neighborhoods.} Indeed, the image of $e^{k+1}$ is
contained in a tubular neighborhood of the image of $e^k$
for each $k\ge1$. In the above construction we consider the
tubular neighborhood of $e^k$ in $\RR^{k+3}$. However since
we proceed inductively using orthogonal bundles of the
successive images of $e^{k+1}, e^{k+2},\dots$, and we
contract the diameter of the tubular neighborhood at each
step by a constant factor $0<\lambda<1$, we have in fact
that \emph{the image of $e^{k+l}$ is in a tubular
  neighborhood of $e^{k}$ for all $k,l\ge1$.}

Therefore, the distance between the image $e^{k+l}$ and the
image of $e^k$ is bounded by $d\sum_{i=k}^{k+l}\lambda^i<d$,
always \emph{inside} a tubular neighborhood $U^{k+l}_k$ of
$e^k(\TT^k)$ in $\RR^{k+l+2}$. Hence there exists a
projection $\pi_k^{k+l}:U^{k+l}_k\to e^k(\TT^k)$ associated
to this tubular neighborhood, for each $k\ge1$ and $l\ge0$.

\subsubsection{The embedding of $\cT^k$ into $B^{k+2}$}
\label{sec:embedd-ct-into}

The previous discussion provides an embedding $e^k$ of
$\TT^k$ into $B^{k+2}$ for each $k\ge1$. Therefore
considering a tubular neighborhood $U^k$ of the compact
submanifold $e^k(\TT^k)$ in $B^{k+2}$ we obtain a projection
$\pi:U\to e(\TT^k)$ such that $\pi^{-1}(w)$ is a
$2$-disk. Thus we obtain an embedding $\hat e^k$ of
$\TT^k\times\DD$ into $B^{k+2}$.

We can assume that the tubular neighborhood has radius
smaller than $\lambda^k d$ and that $U^k\subset U^k_l$ for
all $l<k$.  Then we can consider the projections
$\pi^k:U^k\to e^k(\TT^k)$ and $\pi^k_l:U^k_l\to e^l(\TT^l)$.
In what follows we assume without loss of generality that
$U^k$ is the image of $\hat e^k$.

\subsubsection{The solenoid map through an isotopy of the identity}
\label{sec:soleno-map-through}

The previous construction of the embeddings $e^k$ of
$\TT^k$ depends on the initial embedding $e^1$ of $\TT^1$ on
$B^3$. Moreover it is clear that each of the embeddings
$e^k$ for $k>1$ depend smoothly on $e^1$. Hence a smooth
family $e^1_t$ of embeddings, for $t\in[0,1]$, defines a
smooth family $e^k_t$ of corresponding higher dimensional
embeddings.

\begin{figure}[htpb]
  \centering
\psfrag{e0}{$e^1_0(\TT^1)$}
\psfrag{e1}{$e^1_t(\TT^1)$}
\psfrag{e2}{$e^1_1(\TT^1)$}
\psfrag{U}{$U^1_1$}
  \includegraphics[scale=0.3  ]{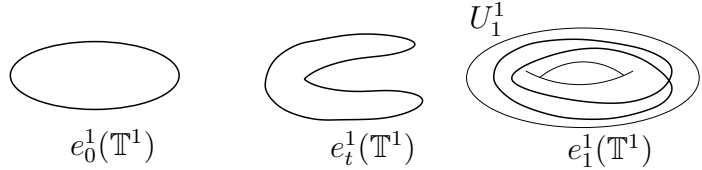}
  \caption{\label{fig:isotopy} The isotopy $e^1_t$ and the
    tubular neighborhood $U^1_1$ in the end.}
\end{figure}

We argue again by induction on $k\ge1$. For $k=1$ we
consider the family of embeddings $e^1_t$ described in
Figure~\ref{fig:isotopy}.

We can construct this family as depicted so that the extreme
elements of the family satisfy $\pi^1_1\circ
e_1^1(\theta)=e_0^1(\theta^2)$ for $\theta\in\TT^1$.  If we
choose a small tubular neighborhood of the image of $e^1_t$,
then we obtain a family of embeddings $\hat e^1_t$ of
$\cT^1$ such that $(\pi^1_1\circ\hat e^1_1)(\theta,z)=
e^1_0(\theta^2)$.

Now we can use the family $e^1_t$ to construct families
$e^k_t$ of embeddings following the inductive procedure
explained before, for each $t\in[0,1]$ and each fixed
$k\ge1$.

Again taking a small tubular neighborhood of the image of
$e^k_t$ we obtain a family of embeddings $\hat e^k_t$ of
$\cT^k$ such that the image of $\hat e^k_1$ of $\cT^k$ is
compactly contained inside the image of $\hat e^k_0$.

Due to the nested construction, for each $k\ge1$ we have $
(\pi^k_1\circ\hat e^k)(\theta_1,\dots,\theta_k,z)=
e^1_0(\theta_1^2) $ after projecting into the lower
dimensional image. Moreover projecting on the previous stage
of the construction we get $ (\pi^k_{k-1}\circ\hat
e^k)(\theta_1,\dots,\theta_k,z)=
e^{k-1}_0(\theta_2^2,\dots,\theta_{k}^2)
\cap(\pi^k_1)^{-1}\{e^1_0(\theta_1^2)\}$, which can easily
be proved by induction on $k\ge1$ following the nested
construction presented above.

This is enough to prove that (independently of the
definition of $A_\Theta$ in $s$)
\begin{align}\label{eq:round}
  \pi^k_k\circ \hat e^k_1 = \pi^k_k \circ \hat e^k_0\circ s
  \quad\text{for all}\quad k\ge1.
\end{align}
Finally, since by definition $\hat e^k_1$ is a small tubular
neighborhood of the image of $e^k_1$, and this set is
contained inside a tubular neighborhood of the image of
$e^k_0$, then from \eqref{eq:round} we see that in fact
there exists a family of contractions
$(A_\Theta)_{\Theta\in\TT^k}$ such that $\hat e^k_1=\hat
e^k_0\circ s$. This completes the proof of
Proposition~\ref{pr:embedding}.

%%%%%%%%%%%%%%%%%%%%%%%%%%%%%%%%%%%%%%%%%%%%%%%%%%%%%

%%% references %%%%%%%%%%%%%%%

\def\cprime{$'$}

% \bibliographystyle{abbrv}
% \bibliography{../bibliobase/bibliography}

\end{document}